\documentclass[11pt,a4paper,reqno]{amsart} 
\usepackage{calc,amsfonts,amsthm,amscd,epsfig,psfrag,amsmath,amssymb,enumerate,graphicx}
\usepackage[showlabels,sections,floats,textmath,displaymath]{}
\reversemarginpar
\newlength\fullwidth
\numberwithin{equation}{section}

\DeclareMathSymbol{\leqslant}{\mathalpha}{AMSa}{"36} 
\DeclareMathSymbol{\geqslant}{\mathalpha}{AMSa}{"3E} 
\DeclareMathSymbol{\eset}{\mathalpha}{AMSb}{"3F}     
\renewcommand{\leq}{\;\leqslant\;}                   
\renewcommand{\geq}{\;\geqslant\;}                   
\newcommand{\liminftwo}[2]{\liminf_{\substack{#1 \\ #2}}} 
\newcommand{\limsuptwo}[2]{\limsup_{\substack{#1 \\ #2}}} 

\def\1{\ifmmode {1\hskip -3pt \rm{I}} \else {\hbox {$1\hskip -3pt \rm{I}$}}\fi}

\newcommand{\var}{\operatorname{Var}} 
\newcommand{\diam}{\operatorname{Diam}} 
 
\newcommand{\la}{\label } \newcommand{\si}{\sigma } 
\newcommand{\be}{\begin{equation} } \newcommand{\tmix}{T_{\rm mix}} 

\newtheorem{Theorem}{Theorem}[section] 
\newtheorem{Lemma}[Theorem]{Lemma} 
\newtheorem{Proposition}[Theorem]{Proposition}

\newtheorem{claim}[Theorem]{Claim}

\newcommand{\N}{\mathbb N}
\newcommand{\bP}{{\bf P}} 
\newcommand{\bE}{{\bf E}}

\newcommand{\cE}{\ensuremath{\mathcal E}} 
 
\newcommand{\cG}{\ensuremath{\mathcal G}}

\newcommand{\cL}{\ensuremath{\mathcal L}}

\newcommand{\cS}{\ensuremath{\mathcal S}}

\newcommand{\bbE}{{\ensuremath{\mathbb E}} }

\newcommand{\bbN}{{\ensuremath{\mathbb N}} } 
 
\newcommand{\bbP}{{\ensuremath{\mathbb P}} } 
 
\newcommand{\bbR}{{\ensuremath{\mathbb R}} }

\newcommand{\bbZ}{{\ensuremath{\mathbb Z}} } 

    \let\d=\delta  \let\e=\varepsilon
 \let\g=\gamma \let\h=\eta    \let\k=\kappa  \let\l=\lambda
      \let\o=\omega      
\let\r=\rho   \let\t=\tau   
  
\let\D=\Delta     \let\L=\Lambda 
\let\O=\Omega      

\def\\{\hfill\break}

\def\thsp{\thinspace}

\def\tthsp{\kern .083333 em}

\def\?{\mskip -10mu}
\newfam\msafam
\newfam\msbfam
\newfam\eufmfam
\def\hexnumber#1{%
  \ifcase#1 0\or 1\or 2\or 3\or 4\or 5\or 6\or 7\or 8\or
  9\or A\or B\or C\or D\or E\or F\fi}
\def\({\left(}
\def\){\right)}

\let\neper=e
\let\ii=i

\def\ie{\hbox{\it i.e.\ }}
\def\eg{\hbox{\it e.g.\ }}

\def\nep#1{ \neper^{#1}}

\def\tc{\thsp | \thsp}
\let\<=\langle
\let\>=\rangle

\def\gap{\mathop{\rm gap}\nolimits}

\newcommand{\wt}{\widetilde }

\begin{document}

\title[]{On the approach to equilibrium for a polymer with adsorption
  and repulsion}

\author{Pietro Caputo}
\address{Dipartimento di Matematica, Universit\`a Roma Tre,
Largo S.\ Murialdo 1, 00146 Roma, Italia. e--mail: {\tt caputo@mat.uniroma3.it}}
\author {Fabio Martinelli}
\address{Dipartimento di Matematica, Universit\`a Roma Tre,
Largo S.\ Murialdo 1, 00146 Roma, Italia. e--mail: {\tt martin@mat.uniroma3.it}}
\author {Fabio Lucio Toninelli}
\address{Ecole Normale Sup\'erieure de Lyon, Laboratoire de Physique
and CNRS, UMR 5672\\ 46 All\'ee
  d'Italie, 69364 Lyon Cedex 07, France.
e--mail: {\tt fltonine@ens-lyon.fr}}

\begin{abstract}
We consider paths of a one--dimensional simple random walk conditioned
to come back to the origin after $L$ steps, $L\in 2\bbN$. In the 
{\em pinning model} each path $\eta$ has a weight $\l^{N(\eta)}$, where
$\l>0$ and $N(\eta)$ is the number of zeros in $\eta$.  When the paths
are constrained to be non--negative, the polymer is said to satisfy a
hard--wall constraint.  Such models are well known to undergo a
localization/delocalization transition as the pinning strength $\l$ is
varied. In this paper we study a natural ``spin flip'' dynamics
for 
these models and derive several estimates on its
spectral gap and mixing time.  In particular, for the system with the
wall we prove that relaxation to equilibrium is always at least as
fast as in the free case (\ie $\l=1$ without the wall), where the
gap and the mixing time are known to scale as $L^{-2}$ and $L^2\log
L$, respectively. This improves considerably over previously known
results.
For the system without the wall we
show that the equilibrium phase transition has a clear dynamical
manifestation: for $\l\geq 1$ relaxation is again at least as fast as
the diffusive free case, but in the strictly delocalized phase
($\l<1$) the gap is shown to be $O(L^{-5/2})$, up to logarithmic
corrections.  As an application of our bounds, we prove stretched
exponential relaxation of local functions in the localized regime.
\\
\textit{Date of submission: September 14, 2007}
\\
\textit{Final version accepted by EJP: February 11, 2008}
\end{abstract}

\maketitle

\thispagestyle{empty}

\section{Introduction}
Consider simple random walk paths 
on $\bbZ$ which start at $0$ and 
end at $0$ after $L$ steps, where $L$ is an even integer, \ie
elements of
$$
\O_L =
\{\eta\in\bbZ^{L+1}:\;
\eta_0=\eta_L=0\,,\;|\eta_{x+1}-\eta_x|=1\,,\;x=0,\dots,L-1\}\,.$$
A well known polymer model (the {\sl pinning model}) is
obtained by assigning to each path $\eta\in\O_L$ a weight
\be
\l^{N(\eta)}\,,
\la{w1}
\end{equation}
where $\l>0$ is a parameter and $N(\eta)$ stands for the number of
$x\in\{1,\dots,L-1\}$ such that $\eta_x=0$, \ie the number of
{\em pinned} sites. If $\l>1$ the weight (\ref{w1}) favors pinning of the path
whereas if $\l<1$ pinning is penalized. The case $\l=1$ is 
referred to as the free case.
Normalizing the weights (\ref{w1}) one has 
a probability measure $\mu=\mu^\l_L$ on the set $\O_{L}$ of all
$\binom{L}{L/2}$ paths. This defines our first polymer model.

The second 
model is obtained by considering only paths that
stay non--negative, \ie elements of $$
\O^+_L =
\{\eta\in\O_L:\;
\eta_x\geq 0\,,\; x=1,\dots,L-1\}\,.$$
Normalizing the weights (\ref{w1}) one obtains  
a probability measure $\mu^+=\mu^{+,\l}_L$ on the set $\O^+_{L}$ of all
$\frac2{L+2}\binom{L}{L/2}$ non--negative paths. 
The positivity
constraint will be often referred to as the presence of a {\em wall}.

\begin{figure}[h]
\centerline{\psfig{file=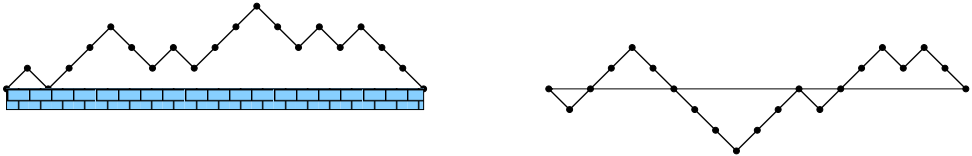,height=0.8in
}}
\caption{Paths with and without the wall, for $L=20$.}
\label{fig:paths_1}
\end{figure}

The two models introduced above have been studied for several decades
and very precise information is available on their asymptotic
properties as $L$ becomes large. The reader is referred to the recent review 
\cite{cf:GB} and references therein and to Section \ref{setup}
below for more details. 
For the moment let us briefly recall that both models display a
transition from a {\em delocalized} to a {\em localized} phase as $\l$ is
increased. Namely, the following  scenario holds. For the 
system without
the wall, if $\l\leq 1$ paths are delocalized (as in the free
case $\l=1$) with $|\eta_{L/2}|$
typically of order $\sqrt{L}$ and a vanishing density of pinned sites, 
while as soon as $\l>1$ paths are strongly localized 
with $|\eta_{L/2}|$ typically of order one with a positive density of
pinned sites. The critical exponents
of the transition can be computed, and the transition itself turns out to 
be of second order: the fraction of pinned sites goes to zero smoothly 
when $\l\searrow 1$.
The system with the wall has a similar behavior but 
the critical point is $\l=2$ instead of $\l=1$. Namely, due to the
entropic repulsion induced by the wall, a small 
reward for pinning (as in the case $1<\l\leq 2$) is
not sufficient to localize the path.

\smallskip

These models and generalizations thereof, where the simple-random-walk
paths are replaced by trajectories of more general Markov chains, are
popular tools in the (bio)-physical literature to describe, e.g.,
pinning of polymers on defect lines in different dimensions, the
Poland-Scheraga model of DNA denaturation, wetting models,...(we refer
for instance to
\cite{cf:Fisher}, \cite[Chap. 1]{cf:GB} and references therein).

Presently there is much activity on the {\sl quenched disordered}
version of these models, where the pinning parameter $\l$ is replaced
by a sequence of (usually log-normal) IID random variables $\l_x,
0<x<L$. The localization-delocalization transition is present also in
this case, and typical questions concern the effect of disorder on the
critical point and on the critical exponents (cf.\ \cite{cf:DHV},
\cite{cf:GT_cmp},
\cite{cf:ken} and \cite{cf:T_qrc}).
Another natural generalization of the polymer models we introduced is
to consider $(d+1)$-dimensional interfaces $\{\eta_x\}_{\{x\in
V\subset \mathbb Z^d\}}$, with or without the hard wall condition
$\{\h_x\geq0\;\forall x\in V\}$, and with some pinning interaction (see
the recent review \cite{cf:velenik} and references therein). 

\bigskip

\begin{figure}[h]
\centerline{
\psfrag{a}{$\frac\l{1+\l}$}
\psfrag{b}{$\frac12$}
\psfrag{0}{$0$}
\psfrag{L}{$L$}
\psfig{file=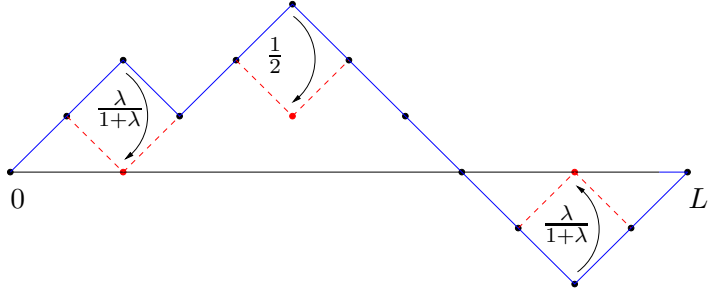,height=1.5in}}
\caption{Three possible transitions, with the corresponding rates, for
the model without the wall.}
\label{fig:dyn}
\end{figure}

We now go back to the two models introduced at the beginning of this
section.  We are interested in the asymptotic behavior of a continuous
time Markov chain naturally associated with them
(cf.\ Figure \ref{fig:dyn}). In the first model
-- system without the wall -- the process is described as follows.
Independently, each site $x\in\{1,\dots,L-1\}$ waits an exponential
time with mean one 
after which the variable $\eta_x$ is updated with
the following rules:
\begin{itemize}
\item if $\eta_{x-1}\neq \eta_{x+1}$, do nothing; 
\item if 
$\eta_{x-1}= \eta_{x+1}=j$ and $|j|\neq 1$, set $\eta_x = j\pm 1$
with equal probabilities; 
\item if $\eta_{x-1}= \eta_{x+1}= 1$, set $\eta_x =0$ with
probability $\frac{\l}{\l+1}$ and $\eta_x =2$ otherwise; 
\item if $\eta_{x-1}= \eta_{x+1}= -1$, set $\eta_x =0$ with
probability $\frac{\l}{\l+1}$ and $\eta_x =-2$ otherwise. 
\end{itemize}

This defines an irreducible Markov
chain on $\O_L$ with reversible probability $\mu$. 
For the system with the wall the process is defined in the same way
with the only difference that now if $\eta_{x-1}= \eta_{x+1}=0$ we are
forced to keep the value $\eta_x=1$. This gives an irreducible Markov
chain on $\O^+_L$ with reversible probability $\mu^+$. 

We shall 
study the speed at which the equilibria $\mu$ and $\mu^+$ are approached 
by our Markov chain 
mostly by way of estimates on the {\em spectral gap} and the 
{\em mixing time}. We refer to
Section \ref{setup} below for the precise definitions, and recall here
that the inverse of the spectral gap (also known as relaxation time)
measures convergence in the $L^2$--norm with respect to the equilibrium
measure, while the mixing time measures convergence in total variation norm 
starting from the worst--case initial condition.

While essentially everything is known about the equilibrium
properties of these polymer models, we feel that there is still much
to understand as far as the approach to equilibrium is concerned. 
In particular, one would like to detect the dynamical signature
of the phase transition recalled above. Our work is a first attempt
in this direction. Before going to a description of our results, we
discuss some earlier contributions.

The problem is well understood in the free case $\l=1$. 
In particular, for the system without the wall, the free case is
equivalent to the so--called symmetric simple exclusion process which
has been analyzed by several authors. We refer to 
the work of Wilson \cite{cf:Wilson}, 
where among other
things the spectral gap of the chain is computed exactly as 
\be
\k_L = 1 - \cos\left(\frac{\pi}L\right)\,,
\la{gap1}
\end{equation}
the principal eigenvalue of the discrete Laplace operator with
Dirichlet boundary conditions, and the mixing time $T_{\rm mix}$ is 
shown to be of order $L^2\log L$ (with
upper and lower bounds differing only by a factor $2$ in the large $L$
limit).

As far as we know, \cite{MaRa_2,MaRa} by Martin and Randall are the
only works where the dynamical problem for all $\l>0$ was
considered. They showed that there is always a polynomial upper bound
on the mixing time. Although their proof is carried out in the case of
the system with the wall only, their result should apply in the
absence of the wall as well. As noted in \cite{MaRa} and as we shall
see in detail in the forthcoming sections, for the system with the
wall, Wilson's coupling method can be easily modified to prove an
upper bound of order $L^2\log L$ on the mixing time for all $\l\leq
1$. On the other hand the problem is harder when $\l>1$, and the
Markov chain decomposition method of \cite{MaRa} only gives $T_{\rm
mix}= O(L^k)$ for some large non--optimal power $k$.

Let us also mention that, on the non-rigorous or numerical level,
various works were devoted recently to the dynamics of polymer models
related to the ones we are considering
(cf.\ for instance 
\cite{cf:dna1,cf:dna2} and references therein). These works are mainly
motivated by the study of the dynamics of heterogeneous DNA molecules
close to the denaturation transition, and therefore focus mainly on
the quenched disordered situation. While the dynamics considered there
is quite different from the one we study here (and in this sense the
results cannot be naturally compared), let us point out that in
\cite{cf:dna1} interesting dynamical transition phenomena are
predicted to occur close to the equilibrium phase transition, both
for the disordered and for the homogeneous models.

\subsection{ Quick survey of our results} We
refer to Section \ref{results} below for the precise statements.
We start with the system with the wall. 
A first result here is that for all $\l>0$, the spectral gap
is bounded below by the gap (\ref{gap1}) of the free case, \ie $\gap
\geq \k_L \sim \pi^2/2L^{2}$. Also, we
prove that for all $\l>0$ the mixing time satisfies $T_{\rm
mix}=O(L^2\log L)$. Furthermore we can prove that these estimates are
optimal (up to constant factors) in the delocalized phase, \ie we can
exhibit complementary bounds for $\l\leq 2$ on the gap and for $\l<2$
on the mixing time. In the localized phase ($\l>2$) we expect the
relaxation to occur faster than in the free case. However, we prove a
general lower bound on the mixing time giving $T_{\rm mix} =
\O(L^2)$
(we recall that by definition $f(x)=\O(g(x))$ for $x\to\infty$ if 
 $\liminf_{x\to\infty}f(x)/g(x)>0$). 
Concerning the spectral gap we show an upper bound $\gap =
O(L^{-1})$. We conjecture these last two estimates to be of the
correct order but a proof of the complementary bounds remains a
challenging open problem\footnote{After this work was completed we
  were able to prove upper and lower bounds on the spectral gap of
  order $L^{-1}$ at least in the perturbative regime $\l=\O(L^{4})$. 
This is part of further work (in progress) on the dynamical aspects of the localization/delocalization transition} (except for $\l=\infty$, where we can
actually prove that $c_1 L^2\leq T_{\rm mix} \leq c_2 L^2$).

The fact that the mixing time grows in every situation at least like
$L^2$ does not exclude that, starting from a particular configuration,
the dynamics can relax to equilibrium much faster. In the localized phase
we explicitly identify such a configuration and show that the dynamics
started from it relaxes within a time $O(\log L)^3$.

Concerning the system without the wall we can show that for all $\l> 1$
the relaxation is at least as fast as in the free case, \ie $\gap\geq
\k_L$
and $T_{\rm mix}= O(L^2\log L)$. However, for $\l>1$ we believe the
true behavior to be the same as described above for $\l>2$ in the
presence of the wall.  On the other hand, the case $\l<1$ is very
different from the system with the wall. Here we prove that the
spectral gap is no larger than $O(L^{-5/2})$, up to logarithmic
corrections,
establishing a clear
dynamical transition from localized to delocalized phase. Describing
the correct asymptotics of the gap (and of the mixing time) for $\l<1$
remains an open problem, although a heuristic argument (see Section
\ref{sec:miclo}) suggests that the $O(L^{-5/2})$ behavior may well be
the correct one.

\begin{table}[ht]
\label{tbl}
\begin{tabular}{|r|c|c|c|}
\hline
\multicolumn{4}{|c|}{Model}\\
\hline
  parameter &  \begin{tabular}{@{}c@{}}conjectured \\ behavior
\end{tabular} & \begin{tabular}{@{}c@{}}rigorous\\lower bound\end{tabular} & \begin{tabular}{@{}c@{}}rigorous\\upper bound\end{tabular} \\
\hline
\hline
\multicolumn{4}{|c|}{Wall, $\l < 2$} \\
\hline
&&&\\[-11pt]
  spectral gap &  $L^{-2}$ & $L^{-2}$ & $L^{-2}$ \\
  mixing time &  $L^2\log L$ & $L^2\log L$ & $L^2\log L$ \\[8pt]
\hline

\multicolumn{4}{|c|}{Wall, $\l=2$} \\
\hline
&&&\\[-11pt]
  spectral gap &  $L^{-2}$ & $L^{-2}$ & $L^{-2}$ \\
  mixing time &  $L^2\log L$ & $L^2$ & $L^2\log L$ \\[8pt]

\hline

\multicolumn{4}{|c|}{Wall, $\l>2$} \\
\hline
&&&\\[-11pt]
  spectral gap &  $L^{-1}$ & $L^{-2}$ & $L^{-1}$ \\
  mixing time &  $L^2$ & $L^2$ & $L^2\log L$ \\[8pt]

\hline
\hline
\multicolumn{4}{|c|}{No wall, $\l<1$} \\
\hline
&&&\\[-11pt]
spectral gap &  $L^{-5/2}$ &  & $L^{-5/2}(\log L)^8$ \\
mixing time &    & $L^{5/2}(\log L)^{-8}$ &  \\[8pt]

\ifodd 1
\hline
\multicolumn{4}{|c|}{No wall, $\l=1$} \\
\hline
&&&\\[-11pt]
 spectral gap & $ L^{-2}$ & $L^{-2}$ & $L^{-2}$ \\
  mixing time &  $L^2\log L$ & $L^2\log L$ &$ L^2\log L$ \\[8pt]

\hline
\multicolumn{4}{|c|}{No wall, $\l>1$} \\
\hline
&&&\\[-11pt]
 spectral gap &  $L^{-1}$ & $L^{-2}$ & $L^{-1}$ \\
 mixing time &  $L^2 $ & $L^2$& $L^2\log L$ \\[8pt]
\hline
\hline
\multicolumn{4}{|c|}{Wall/No wall, $\l=+\infty$} \\
\hline
&&&\\[-11pt]
  mixing time &  $L^2 $ & $L^2$& $L^2$ \\[8pt]

\fi
\hline
\end{tabular}
\smallskip\smallskip\smallskip\smallskip\smallskip
\caption{Rough summary of spectral gap and mixing time bounds. All the
  entries in the table have to be understood as valid up to
  multiplicative constants
  independent of $L$. The statements of our theorems clarify whether 
the bounds hold with constants depending on
$\l$ or not. Blank entries in the table correspond to questions which
have not been addressed in this work.}
\end{table}
Finally, besides focusing on global quantities like gap and mixing
time, it is of interest to study how local observables, \eg the
local height function $\h_x$, relax to equilibrium. Note that this
point of view is closer to the one of the theoretical physics papers
\cite{cf:dna1,cf:dna2} we mentioned above. This question is particularly
interesting in the localized phase, where the infinite-volume
equilibrium measure 
is the law of a  positive recurrent Markov chain and $\eta_x$ is of order
one. As a consequence of 
the fact that the
spectral gap vanishes for $L\to\infty$ as an inverse power of $L$, we
will show in Theorem \ref{th:t13} upper and lower
bounds of stretched exponential type for the relaxation of local
functions.

\smallskip
The work is organized as follows: in Section \ref{sec:eq} the model is
defined and some basic equilibrium properties are recalled; in Section
\ref{sec:dyn}
we introduce our dynamics and for completeness we define a few
standard tools (spectral gap, mixing time, etc.); in Section
\ref{firstarg} 
we describe a basic coupling argument due to D. Wilson
\cite{cf:Wilson}, which we use at various occasions; in Section
\ref{results} we state our main results, which are then proven in
Sections \ref{up} to \ref{locco}.

\section{Setup and preliminaries}\la{setup}

In this section we set the notation and collect several tools to be
used repeatedly in the rest of the paper.  
  
\subsection{Some equilibrium properties}
\label{sec:eq}
Fix $\l>0$ and $L\in 2\bbN$ and write $\L:=\{0,\dots,L\}$. As in the introduction $\mu=\mu^\l_L$
denotes the equilibrium measure of the unconstrained system. 
The Boltzmann weight associated to a configuration $\eta\in \Omega_L$
is
\begin{eqnarray}
\label{eq:Boltz_nowall}
\mu^{\lambda}_L(\eta):=\frac{\lambda^{N(\eta)}}{Z_L(\lambda)}\,,
\end{eqnarray}
where $N(\eta):=\#\{0<x<L:\eta_x=0\}$ and
\begin{eqnarray}
Z_L(\lambda):=\sum_{\eta\in \Omega_L} \lambda^{N(\eta)}\,.
\end{eqnarray}

The equilibrium of the constrained system is described by
$\mu^+=\mu^{+,\lambda}_L$. Here the Boltzmann weight associated
to a configuration $\eta\in\Omega_L^+$ is
\begin{eqnarray}
\label{eq:Boltz_wall}
\mu^{+,\lambda}_L(\eta):=\frac{\lambda^{N(\eta)}}{Z^+_L(\lambda)}\,,
\end{eqnarray}
where
\begin{eqnarray}
Z^+_L(\lambda):=\sum_{\eta\in\Omega^+_L} \lambda^{N(\eta)}\,.
\end{eqnarray}
When there is no danger of confusion, we will omit the 
indexes $\lambda$ and $L$ and write $\mu$ for $\mu_L^\l$ and $\mu^+$ for $\mu^{+,\lambda}_L$. 

Considering reflections of the path between consecutive zeros one obtains 
the following identity:
\begin{eqnarray}
\label{eq:identita}
2\,  Z^+_L(2\lambda)= Z_L(\lambda)\,.
\end{eqnarray}
Moreover, if $\zeta(\eta):=\{x\in \L:\eta_x=0\}$ is the set of zeros of
the configuration $\eta$, one has 
\begin{eqnarray}
\label{eq:id2}
  \mu^{+,2\lambda}_L(\zeta= S)=  \mu^{\lambda}_L(\zeta=S)\,,\quad S\subset \L\,.
\end{eqnarray}
In other words, the thermodynamic properties of the two models are 
essentially equivalent modulo a change of $\l$. On the other hand, 
we will see that the two present very different dynamical phenomena.

\subsubsection{Free energy and the localization/delocalization transition}
Let $\bP$ and
$\bE$ denote the law and expectation of the one--dimensional simple
random walk $\eta:=\{\eta_n\}_{n\geq 0}$ with initial condition
$\eta_0=0$. Then,
\begin{eqnarray}
 Z_L(\lambda)=2^L \bE \left(\lambda^{N(\eta)}\,1_{\{\eta_L=0\}}\right)\,,
\end{eqnarray}
and 
\begin{eqnarray}
  Z^+_L(\lambda)=2^L \bE \left(\lambda^{N(\eta)}\,1_{\{\eta_L=0\}}
  \,1_{\{\h_x\geq 0\;\forall x<L\}}\right)\,.
\end{eqnarray}
The {\sl free energy} is defined for the system without the wall
as 
\begin{eqnarray}
\label{eq:F}
F(\l):=\lim_{L\to\infty}\frac1L \log Z_L(\l)-\log 2\,.
\end{eqnarray}
The limit exists by super-additivity. Similarly, the free energy of the 
system with the wall is denoted by $F^+(\l)$. 
Of course, one has $F^+(\l)= F(\l/2)$, as follows from
\eqref{eq:identita}.

The following is well known (cf.\ \eg \cite[Ch. 2]{cf:GB}): 
$F(\l)=0$ for $\l\leq 1$ and $F(\l)>0$ for $\l>1$. Moreover, for
$\l>1$, $F(\l)$ can be equivalently defined as the unique
positive solution of
\begin{eqnarray}
\label{eq:Faltern}
  \sum_{n\in 2\N}\bP(\inf\{k>0:\eta_k=0\}=n)\nep{-n\,F}=\frac1\lambda.
\end{eqnarray}
Together with the explicit expression for the Laplace transform of the
first return time of the simple random walk,
\begin{eqnarray}
\sum_{n\in 2\N} z^n \bP(\inf\{k>0:\eta_k=0\}=n)=1-\sqrt{1-z^2}
\end{eqnarray}
for $|z|\leq 1$,  \eqref{eq:Faltern} implies 
\begin{eqnarray}
\label{eq:tildeFl1}
F(\l)=\frac12 
\log\left[\frac{\l^2}{2\l-1}\right]\,,  
\end{eqnarray}
 for $\l>1$. Note that $F^+(\l)>0$ if and only if $\l>2$.

We will need the following sharp estimates on the asymptotic 
behavior of the partition function for large $L$:
\begin{Theorem}\cite[Th. 2.2]{cf:GB}
\label{th:asymptZ}
\begin{eqnarray}
\label{eq:asymptZ}
  2^{-L}Z_L(\lambda)
 \stackrel{L\to\infty}\sim C(\lambda) \times
\left\{
\begin{array}{lll}
\nep{L\, F(\lambda)} & \mbox{for} & \lambda>1\\
 L^{-1/2} & \mbox{for} & \lambda=1\\
 L^{-3/2} & \mbox{for} & \lambda<1
\end{array}
\right.
\end{eqnarray}
where $C(\lambda)>0$ for every $\lambda$, \ie the ratio of the two
sides in (\ref{eq:asymptZ}) converges to one.
\end{Theorem}
We refer to \cite[Th. 2.2]{cf:GB} for an expression of $C(\lambda)$ in terms
of the law $\bP(\cdot)$. 
From the explicit expression \eqref{eq:tildeFl1} one sees that
$F(\cdot)$ is differentiable with respect to $\l$ in
$(0,\infty)$. Since the free energy is a convex function of $\log \l$,
one deduces that the average density of pinned sites satisfies
\begin{eqnarray}
\label{eq:density}
\lim_{L\to\infty} \frac1L \mu^\l_L(N(\eta))=\frac{d F(\l)}
{d\log \l}
\;\left\{
\begin{array}{lll}
=0&\mbox{if}&\l\leq 1\\
>0&\mbox{if}&\l>1
\end{array}
\right .
\end{eqnarray}
For this reason, one calls  the region of parameters $\l\leq1$  {\sl delocalized
phase} and $\l>1$ {\sl localized phase}, and $\l=1$ the {\sl critical point}
(for the system with the wall, the critical point is therefore $\l=2$).

One can go much beyond the density statement \eqref{eq:density}
in characterizing the two phases. In the rest of this section we recall
some known results.

\subsubsection{The strictly delocalized phase} 
\label{sec:delocph}
This terminology refers
to the situation $\l<1$ (or $\l<2$ with the wall). In this, case, the number
of zeros $N(\eta)$ is typically finite and its law has an exponential tail.
In what follows we write $c=c(\l)$ for a positive 
constant (not necessarily the same at each occurrence)
which can depend on $\l$ but not on $L$.
There exists $c=c(\l)$ such that 
\begin{eqnarray}
\label{eq:pochizeri}
\mu^\l_L(N(\eta)\geq j)\leq c\, \nep{- j/c}\,,
\end{eqnarray}
uniformly in $L$. (This simply follows from
\begin{eqnarray}
\mu^\l_L(N(\eta)\geq j)\leq e^{-\varepsilon j}\mu^\l_L\left(
e^{\varepsilon N(\eta)}\right)=e^{-\varepsilon j}\frac{Z_L(\l\, e^\varepsilon)}
{Z_L(\l)}\,,
\end{eqnarray}
if we choose $\varepsilon>0$ small enough so that $\l\exp(\varepsilon)<1$,
cf.\ Theorem \ref{th:asymptZ}.)
It is also easy to see that there is a non-zero probability that $N(
\eta)=0$:
\begin{eqnarray}
\label{eq:nozeri}
\mu_L^\l(N(\h)=0)=2\frac{
\bP(\eta_L=0,\;\eta_x>0\;\;\forall\; 1<x<L)}{2^{-L}
Z_L^\l}\stackrel{L\to\infty}\sim c\in (0,1)\,,
\end{eqnarray}
where in the last step we used \eqref{eq:asymptZ} and the fact that
\begin{eqnarray}
\label{eq:denomina}
\lim_{L\to\infty}L^{3/2}\bP(\eta_L=0,\;\eta_x>0\;\;\forall\; 1<x<L)>0\,,
\end{eqnarray}
\cite[Sec. III.3]{cf:Feller}.
Finally, we will need the following upper bound on the probability
that there exists a zero far away from the boundaries of the system:
\begin{eqnarray}
\label{eq:farfrom0L}
\mu _L^\l
(\exists x:\, \ell\leq  x\leq L-\ell,\eta_x=0)\leq \frac{c}{\ell^{1/2}},
\end{eqnarray}
for every $L$ and $\ell<L/2$.
This can be extracted immediately from Theorem \ref{th:asymptZ}.

\subsubsection{The localized phase}
Here $\l>1$ for the system without the wall or $\l>2$ with the wall.
In the localized phase, $|\h_x|$ is typically of order $1$ with
exponential tails, and correlation functions between local functions
decay exponentially fast. Given a function $f:\Omega_L\rightarrow\bbR$
we denote by $\mathcal S_f$ the support of $f$, \ie the minimal set
$I\in \L$ such that $f$ depends only on $\{\h_x\}_{x\in I}$, and set
$\|f\|_\infty:=\max_{\eta\in \Omega_L} |f(\eta)|$.  Then, it is not
difficult to prove:
\begin{Lemma}
\label{th:tight}
Let $\l>1$. For every $L\in2\bbN$ and $x,\ell\leq L$
\begin{eqnarray}
\label{eq:tight}
\mu_L^\l(|\eta_x|\geq \ell)\leq c\, \nep{- \ell F(\l)}.
\end{eqnarray}
Moreover, for every pair of  functions $f,g:\Omega_L\to\bbR$
\begin{eqnarray}
\label{eq:decay}
\left|\mu^\l_L(f\,g)-\mu^\l_L(f) \mu_L^\l(g)
\right|\leq c\, \|f\|_\infty\,\|g\|_\infty \nep{-d(\mathcal S_f,\mathcal S_g)/c}
\end{eqnarray}
where $d(\cdot,\cdot)$ denotes the usual distance between subsets of
$\bbZ$.  One has exponential loss of memory of boundary conditions:
\begin{eqnarray}
\label{eq:bc}
\sup_{L>k} \left|\mu^\l_L(f)-\mu^\l_k(f)\right|\leq c\, \|f\|_\infty
\nep{- d(\mathcal S_f,\{k\})/c},
\end{eqnarray}
where $d(\mathcal S_f,\{k\})$ is the distance between 
$\cS_f\subset \{0,\dots,k\}$ and the point
$\{k\}$. 
Finally, for every bounded local function the thermodynamic limit
\begin{eqnarray}
\lim_{L\to\infty}\mu_L^\l(f)
\end{eqnarray}
exists. The same holds for $\mu^{+,\l}_L$ if $\l>2$.
\end{Lemma}
These results follow for instance from those proven in
\cite{cf:GT_alea} in a more general context, \ie when the constant
$\l$ is replaced by a sequence of IID random variables $\l_x,x\in\L$.

\subsection{The Markov chain}
\label{sec:dyn}
The process described in the introduction is nothing but the standard
heat bath dynamics. For the system without the wall we can formulate
this as follows. Let $Q_x$ denote the $\mu$--conditional expectation at
$x$ given the values of the heights $\eta_y$ at all vertices $y\neq
x$, where $\mu=\mu^\l_L$ is the equilibrium measure (\ref{eq:Boltz_nowall}). 
Namely, for all $f:\O_L\to\bbR$, and $x\in\{1,\dots,L-1\}$ we write
\be
Q_xf = \mu(f\tc \eta_y\,,\;y\neq
x)\,.
\la{qx}
\end{equation}
Our process is then the continuous-time Markov chain with
infinitesimal
generator given by 
\be
\cL f = \sum_{x=1}^{L-1} \left[Q_xf - f\right]\,,\;\quad f:\O_L\to\bbR\,.
\la{gene}
\end{equation}
Note that the generator can be written in more explicit terms as 
$$
\cL f(\eta) = \sum_{x=1}^{L-1} c_x(\eta)\left[f(\eta^x) - f(\eta)\right]\,,
$$
where $\eta^x$ denotes the configuration $\eta$ after the $x$-th
coordinate has been ``flipped'', and the rates $c_x(\eta)$ are given
by
$$
c_x(\eta)=\begin{cases}
\frac12 & \;\eta_{x-1}=\eta_{x+1}\notin\{-1,1\}\\
\frac\l{\l+1} & \;(\eta_{x-1},\eta_x,\eta_{x+1}) = (1,2,1) \;\text{or}\;
(-1,-2,-1)\\
\frac1{\l+1} & \;(\eta_{x-1},\eta_x,\eta_{x+1}) = (1,0,1) \;\text{or}\;
(-1,0,-1)\\
0 & \; \eta_{x-1}\neq\eta_{x+1}
\end{cases}
$$ 
We shall write $P_t$, $t\geq 0$, for the associated semigroup acting
on functions on $\O_L$. Given an initial condition $\xi$, we write 
$\eta^\xi(t)$ for the configuration at time $t$, so that the expected
value of $f(\eta^\xi(t))$ can be written as $P_tf (\xi)$. 
 
Similarly, in the presence of the wall, 
if $Q_x^+$ denotes the $\mu^+$--conditional expectation at
$x$ given the path at all vertices $y$, $y\neq
x$, where $\mu^+=\mu^{+,\l}_L$ is the equilibrium measure (\ref{eq:Boltz_wall}), the infinitesimal generator becomes
\be
\cL^+ f = \sum_{x=1}^{L-1} \left[Q^+_xf - f\right]\,,\;\quad f:\O_L^+\to\bbR\,.
\la{gene+}
\end{equation}
We write  $\eta^{+,\xi}(t)$ for  the configuration at time $t$ with
initial condition $\xi$. 
Similarly, we write $P_t^+$ for the associated semigroups acting
on functions on $\O_L^+$. 
If no confusion arises we shall drop the $+$ superscript and
use again the notation $\eta^\xi(t)$, $P_t$ as in the case without the wall.

\subsubsection{Coupling and monotonicity}\la{coupling}
A standard procedure allows to define a probability measure
$\bbP$ which is a simultaneous coupling of the laws of processes associated
to different initial conditions. Moreover, the measure $\bbP$ can be
used to couple the laws of processes corresponding to different values of $\l$
and to couple paths evolving with the wall to paths evolving without
the wall. 

The construction of $\bbP$, the {\em global coupling}, 
can be described as follows. We need $L-1$ independent Poisson
processes $\o_x$ with parameter $1$, which mark the updating times at each
$x\in\{1,\dots,L-1\}$, and a sequence $\{ u_n\,,\;n\in\bbN\}$
of independent random variables with uniform distribution in
$[0,1]$, which stand for the ``coins'' to be flipped for the updating
choices. 
Given an initial condition $\xi$, a realization $\o$ of the 
Poisson processes and a realization $u$ of the variables $u_n$ we can
compute the path $\eta^\xi(s)$, $s\leq t$, for any fixed $t>0$, 
as follows: sites to be updated together with their updating
times up to time $t$ are chosen according to $\o$; if 
the $k$-th update occurs at site $x$ and at time $s_k$, and
$\eta^\xi_{x-1}(s_k)=\eta^\xi_{x+1}(s_k)=j$ 
then
\begin{itemize}
\item if $|j|\neq 1$, set $\eta_x=j+1$ if $u_k\leq \frac12$,
and $\eta_x=j-1$ otherwise; 
\item if $j= 1$, set $\eta_x=0$ if $u_k\leq \frac\l{\l+1}$,
and $\eta_x=2$   otherwise; 
\item if $j= -1$, set $\eta_x=0$ if $u_k\leq \frac\l{\l+1}$,
and $\eta_x=-2$   otherwise.
\end{itemize}
Of course, in case of an evolution with the wall we have to add the 
constraint that a site $x$ such that
$\eta^\xi_{x-1}(s_k)=\eta^\xi_{x+1}(s_k)=0$ cannot change.

We can run this process for any initial data $\xi$. 
It is standard to check that, provided we use the {\em same}
realization $(\o,u)$ for all copies, the above construction 
produces the desired coupling.

Given two paths
$\xi,\si\in\O_L$ we say that $\xi\leq \si$ iff $\xi_x\leq \si_x$
for all $x\in\L$. 
By construction, if $\xi\leq \si$, 
then $\bbP$--a.s.\ 
we must have $\eta^\xi(t)\leq \eta^\si(t)$ at all times. The same
holds for the evolution with the wall. 
In particular, we will be interested in the evolution started from the
{\em maximal} path $\wedge$, defined as $\wedge_x=x$ for $x\leq L/2$
and $\wedge_x=L-x$ for $L/2\leq x\leq L$, and from the {\em minimal}
path $\vee:=-\wedge$. For the system with the wall the
minimal path is the zigzag line given by $\eta_x=0$ for all even $x$ 
and $\eta_x=1$ for all odd $x$. For simplicity, we shall again use the
notation $\vee$ for this path.

Note that if the initial condition $\xi$ evolves with the wall
while $\si$ evolves without the wall we have $\eta^\si(t)\leq
\eta^{+,\xi}(t)$, if $\si\leq \xi$. Finally, for evolutions with the
wall we have an additional monotonicity in $\l$, \ie if 
$\si$ evolves with parameter $\l$ and $\xi $ with parameter $\l'$
then $\eta^{+,\si}(t)\leq \eta^{+,\xi}(t)$ if $\si\leq \xi$ and
$\l\geq \l'$. 

Let $\bbE$ denote expectation with respect to the global coupling
$\bbP$. Using the notation $\bbE[f(\eta^\xi(t))] = P_t f(\xi)$ the
monotonicity discussed in the previous paragraph takes the form of the
statement that for every fixed $t\geq 0$, the function $P_t f$ is
increasing whenever $f$ is increasing, where a function $f$ is called
{\em increasing} if $f(\xi)\geq f(\si)$ for any $\si,\xi$ such that
$\si\leq \xi$. A whole family of so--called FKG inequalities can be
derived from the global coupling. For instance, the comparison between
different $\l$'s mentioned above, by taking the limit $t\to\infty$
yields the inequality $\mu^{+,\l}(f)\leq
\mu^{+,\l'}(f)$, valid for any increasing $f$ and any $\l\geq \l'$. 
Also, a straightforward modification of the same argument 
proves that for any subset $S\subset \L$ and
any pair of paths
$\si,\xi\in\O_L$ such that $\si\leq \xi$, then 
\be
\mu(f\tc \eta=\si \;\text{on}\;S)\leq \mu(f\tc \eta=\xi
\;\text{on}\;S)\,,
\la{mono1}
\end{equation}
for every increasing $f:\O_L\to\bbR$. The same arguments apply in the
presence of the wall, giving (\ref{mono1}) with $\mu^+$ in place of
$\mu$,
for every increasing $f: \O_L^+\to\bbR$. 

We would like to stress that monotonicity and its consequences such as
FKG inequalities play an essential role 
in the analysis of our models. Unfortunately, these nice properties need not be
available in other natural polymer models.

\subsubsection{Spectral gap and mixing time}
To avoid repetitions we shall state the following definitions for the system 
without the wall only (otherwise simply replace $\mu$ by $\mu^+$,
$\cL$ by $\cL^+$ etc.\ in the expressions below). 

Let $P_t(\xi,\xi')=\bbP(\eta^\xi(t)=\xi')$ denote the 
kernel of our Markov chain. It is easily checked that $P_t$ satisfies 
reversibility with respect to $\mu$, \ie
\be
\mu(\xi)P_t(\xi,\xi') = \mu(\xi')P_t(\xi',\xi)\,,\quad \xi,\xi'\in\O_L\,,
\la{reve}
\end{equation}
or, in other terms, $\cL$ and $P_t$ are self--adjoint in $L^2(\mu)$.
In particular, $\mu$ is the unique invariant distribution and 
$P_t(\xi,\eta)\to \mu(\eta)$ as $t\to\infty$ for every $\xi,\eta\in\O_L$. 
The rate at which this convergence takes place will be  
measured using the following standard tools.

The Dirichlet form associated to (\ref{gene}) is:
\begin{align}
\label{eq:dirich}
\cE(f,f)=
- \mu(f\cL f) =\sum_{0<x<L}
\mu\left[(Q_xf-f)^2\right]\,.
\end{align}
The spectral gap is defined by
\begin{eqnarray}
\gap=\inf_{f:\O_L\to\bbR}
\frac{\cE(f,f)}{\var(f)}\,,
\end{eqnarray}
where $\var(f)=\mu(f^2)-\mu(f)^2$ denotes the variance.
The spectral gap is the smallest non--zero eigenvalue of $-\cL$. 
It measures the rate of exponential decay
of the variance 
of $P_t f$ as $t\to \infty$, \ie $\gap$ is
the (optimal) constant such that for any $f$, $t>0$:
\be
\var(P_tf)\leq \nep{-2t\gap}\,\var(f)\,.
\la{gapvar}
\end{equation}

The mixing time $\tmix$ 
is defined by
\begin{eqnarray}
\tmix=\inf\{t>0 \,:\;\max_{\xi\in \Omega_L}
\|P_t(\xi,\cdot)-\mu\|_{\rm var}
\leq 1/{e}\}\,,
\end{eqnarray}
where $\|\cdot\|_{\rm var}$ stands for the usual total variation norm:
$$
\|\nu-\nu'\|_{\rm var} =
\frac12\sum_{\eta\in\O_L}|\nu(\eta)-\nu'(\eta)|\,,$$
for arbitrary probabilities $\nu,\nu'$ on $\O_L$. 
We refer \eg to Peres \cite{Peres} for more background on
mixing times.
Using familiar relations between total variation distance and
coupling and using the monotonicity of our Markov chain 
we can estimate, for any $\xi$ and $t>0$:
\be
\|P_t(\xi,\cdot)-\mu\|_{\rm var} \leq 
\bbP\left(\eta^\wedge(t)\neq \eta^\vee(t)\right)\,,
\la{varcoup}
\end{equation}
where $\eta^\wedge(t),\eta^\vee(t)$ denote the evolutions from maximal and
minimal paths respectively. 
This will be our main tool in estimating
$\tmix$ from above. Also, (\ref{varcoup}) will be used to 
estimate the spectral gap from below. Indeed, a standard argument 
(see \eg Proposition 3 in \cite{cf:Wilson}) shows that 
$- \liminf_{t\to\infty} \frac1{t}\log \left(\max_{\xi}
\|P_t(\xi,\cdot)-\mu\|_{\rm var}\right)$ is a lower bound on the gap, so that
\be
\gap \geq - \liminf_{t\to\infty} \frac1{t}\log \bbP\left(\eta^\wedge(t)\neq \eta^\vee(t)\right)\,.
\la{gapcoup}
\end{equation}

Finally, it is well known 
that $\gap$ and $\tmix$ satisfy the general relations
\be
\gap^{-1}\leq \tmix\leq \gap^{-1}(1-\log\mu_*)\,,
\la{gaptmix}
\end{equation}
where $\mu_*=\min_\eta\mu(\eta)$. Note that 
in our case $-\log \mu_*=O(L)$ for every fixed $\l$.

\subsection{A first argument}\la{firstarg}
Let $\D$ denote the 
discrete Laplace operator $$(\D\varphi)_x=\frac12
(\varphi_{x-1}+\varphi_{x+1}) - \varphi_x\,.$$ 
We shall need the following computation in the sequel.
\begin{Lemma}
\la{lemma1}
Set $\d = 2/(1+\l)$. 
For the system without the wall, for every $x=1,\dots,L-1$:
\be
\cL\eta_x = (\D\eta)_x + (1-\d) \,1_{\{\eta_{x-1} = \eta_{x+1} = - 1\}}
- (1-\d) \,1_{\{\eta_{x-1} = \eta_{x+1} =  1\}}\,.
\la{lemma10}
\end{equation}
For the system with the wall, for every $x=1,\dots,L-1$:
\be
\cL^+\eta_x = (\D\eta)_x + 1_{\{\eta_{x-1} = \eta_{x+1} = 0\}}
- (1-\d) \,1_{\{\eta_{x-1} = \eta_{x+1} =  1\}}\,.
\la{lemma11}
\end{equation}
\end{Lemma}
If $\l=1$, then $\d=1$ so that (\ref{lemma10}) has pure diffusive character. 
If $\l\neq 1$ the correction terms represent
the attraction ($\l>1$) or repulsion ($\l<1$) at $0$. In the presence of the wall
there is an extra repulsive term.  
\proof
From (\ref{gene}) we see that $\cL\eta_x = \mu[\eta_x\tc \eta_{x-1},\eta_{x+1}] -\eta_x$. If $\eta_{x-1}\neq \eta_{x+1}$ then 
$\mu[\eta_x\tc \eta_{x-1},\eta_{x+1}] = \frac12(\eta_{x-1}+\eta_{x+1})$.
The same holds if $\eta_{x-1} =\eta_{x+1} = j$ with $|j|\neq 1$. 
Finally, if $\eta_{x-1} =\eta_{x+1} = \pm 1$ we have that 
$$\mu[\eta_x\tc \eta_{x-1},\eta_{x+1}] = \pm \d = \d\,\frac12(\eta_{x-1}+\eta_{x+1})\,.$$ 
This proves (\ref{lemma10}). 
The proof of (\ref{lemma11}) is the same, with the observation that 
$$\mu[\eta_x\tc \eta_{x-1},\eta_{x+1}] = 1\,,$$ if $\eta_{x-1} =\eta_{x+1} = 0$. 
\qed

\smallskip

Next, we describe an argument which is at the heart of Wilson's
successful analysis of the free case $\l=1$. 
Define the non-negative profile function $g_x:=\sin\left(\frac{\pi x}L\right)$ and
observe that $g$ satisfies
\be
(\D g)_x = - \k_L \,g_x\,,\quad x\in\{1,\dots,L-1\}\,,
\la{gDg}
\end{equation}
where $\k_L$ is the first Dirichlet eigenvalue of $\D$ given in
(\ref{gap1}). Define \be
\Phi(\eta)=\sum_{x=1}^{L-1} g_x \eta_x\,.
\la{phi}
\end{equation}
Lemma \ref{lemma1} shows that for $\l=1$, for the system without the wall, one has 
\be
\cL\Phi=\sum_{x=1}^{L-1} g_x (\D\eta)_x = \sum_{x=1}^{L-1} (\D g)_x \eta_x
= -\k_L\Phi\,,
\la{arg30}
\end{equation}
where we use summation by parts and (\ref{gDg}). Therefore 
$P_t\Phi(\eta)=\nep{-\k_L t }\,\Phi(\eta)$ for all $t$ and $\eta$.
In particular, if we define 
\be
\widetilde\Phi_t = \sum_{x=1}^{L-1} g_x (\eta^\wedge_x(t) - \eta^\vee_x(t))\,,
\la{phit}
\end{equation}
then 
$\bbE\widetilde\Phi_t = P_t\Phi(\wedge)-P_t\Phi(\vee) = 
\widetilde\Phi_0\,\nep{-\k_L t }$. Note that monotonicity implies that $\widetilde\Phi_t\geq 0$ for all $t\geq 0$. 
Since $g_x\geq \sin(\pi/L)$, $0<x<L$, we have 
\begin{align}
\bbP\left(\eta^\wedge(t)\neq \eta^\vee(t)\right) &
\leq \bbP\left(\widetilde\Phi_t\geq 2\sin(\pi/L)\right) \nonumber\\
&\leq \frac{\bbE\widetilde\Phi_t}{2\sin(\pi/L)} = \frac{\widetilde\Phi_0\nep{-\k_L t }}{2\sin(\pi/L)}\la{phit1}
\end{align}
Inserting (\ref{phit1}) in (\ref{gapcoup}) one obtains 
\be
\gap\geq \k_L\,.
\la{arg1}
\end{equation}
(Since here $\cL\Phi=-\k_L\Phi$ this actually gives $\gap= \k_L$.)
Using (\ref{varcoup}) 
one has the upper bound 
$\tmix\leq \k_L^{-1}\log\frac{\nep{}\,\widetilde\Phi_0}{2\sin(\pi/L)}$.
Since $\k_L\sim \pi^2/2L^2$ and $\widetilde\Phi_0\leq L^2/2$, we have 
\be
\tmix\leq \left(\frac{6}{\pi^2} + o(1)\right)\,L^2\log L\,.
\la{arg2} 
\end{equation}
The estimate (\ref{arg2}) is of the correct order in $L$, although the constant
might be off by a factor $6$, cf.\ Wilson's work \cite{cf:Wilson} for more details.

\section{Main results}\la{results}

\subsection{Spectral gap and mixing time with the wall}
The first result shows that relaxation will never be slower than in the free
case without the wall. 
\begin{Theorem}\la{th1}
For  every $\l>0$, 
\be
\gap \geq \k_L\,,
\la{th1_1}
\end{equation}
where $\k_L = 1 - \cos\left(\frac{\pi}L\right)$. Moreover,
\be
\tmix\leq \left(\frac{6}{\pi^2} + o(1)\right)\,L^2\log L\,.
\la{th1_2} 
\end{equation}
\end{Theorem}
The proof of these estimates will be based on a comparison with the 
free case, which boils down to a suitable control on the correction terms
described in Lemma \ref{lemma1}. This will be worked out in Section \ref{up}.

The next theorem gives complementary bounds which imply that Theorem \ref{th1}
is sharp up to constants in the strictly delocalized phase.
\begin{Theorem}\la{th2}
For  every $\l\leq 2$, 
\be
\gap \leq c\,L^{-2}\,,
\la{th2_1}
\end{equation}
where $c>0$ is independent of $\l$ and $L$. Moreover, for $\l<2$ we have  
\be
\tmix\geq \left(\frac1{2\pi^2}+ o(1)\right)\,L^2\log L\,.
\la{th2_2} 
\end{equation}
For $\l>2$ we have 
\be
\gap \leq c\,L^{-1}\,,
\la{th2_3}
\end{equation}
where $c=c(\l)$ is independent of $L$. Finally, for every $\l > 0$: 
\be
\tmix\geq c\,L^2\,,
\la{th2_4} 
\end{equation}
for some $c>0$ independent of $\l$ and $L$.
\end{Theorem}

The proof of the upper bounds (\ref{th2_1}) and (\ref{th2_3}) 
will be obtained by choosing a suitable test function in the
variational principle defining the spectral gap. 
The estimate (\ref{th2_2}) will be achieved by a suitable comparison
with the free case, while (\ref{th2_4}) will follow by a comparison with the
extreme case $\l=\infty$. 
These results are proven in Section \ref{down}.

We expect the $L^2\log L$ estimate (\ref{th2_2}) to hold at the critical 
point $\l=2$ as well, but for our proof we require strict delocalization
(in (\ref{th2_2}) what may depend on $\l$ is the $o(1)$ function).

We conjecture the estimates (\ref{th2_3}) and (\ref{th2_4}) to be sharp (up
to constants) in the localized phase $\l>2$. In particular, in
Proposition \ref{lainf} we prove
that (\ref{th2_4}) is sharp at $\l=\infty$. 

It is interesting that, although the mixing time is 
$\Omega(L^2)$ in every situation, for the model with the wall we can prove that
the dynamics converges to the invariant measure much faster if started
from the minimal configuration, $\vee$, which so to speak is already
``sufficiently close to equilibrium'':
\begin{Theorem}
\label{th:relaxV}
For $\l>2$ there exists $c(\l)<\infty$ such that 
\begin{eqnarray}
\label{eq:UB_V}
\limsuptwo{L\to\infty\,,}{t\geq c(\l)(\log L)^3}
\| P_t(\vee,\cdot)-\mu_L^{+,\l}\|_{\rm var}=0.
\end{eqnarray}
On the other hand 
\begin{eqnarray}
\label{eq:LB_V}
\liminftwo{L\to\infty\,,}{t\leq(\log L)^2/c(\l)}
\| P_t(\vee,\cdot)-\mu_L^{+,\l}\|_{\rm var}=1.
\end{eqnarray}
\end{Theorem}

The proof of Theorem \ref{th:relaxV} can be found in Section
\ref{locco}.

\subsection{Spectral gap and mixing time without the wall}
We start with the lower bounds on the gap and upper bounds on $\tmix$.
\begin{Theorem}\la{th4}
For any $\l \geq 1$, $\gap$ and $\tmix$ satisfy (\ref{th1_1}) and (\ref{th1_2}) respectively. 
\end{Theorem}
The proof is somewhat 
similar to the proof of Theorem \ref{th1} 
and it will be given in Section \ref{up}. 
We turn to the upper bounds on the gap and lower bounds on $\tmix$. 
\begin{Theorem}\la{th5}
For $\l>1$, $\gap$ and $\tmix$ satisfy 
(\ref{th2_3}) and (\ref{th2_4}) respectively. If $\l<1$, on the other hand, 
there exists $c(\l)<\infty$ such that
\begin{eqnarray}
  \gap\leq {c(\l)}\frac{(\log L)^8}{L^{5/2}}\,.
  \la{th5_2}
\end{eqnarray}
\end{Theorem}
The proof of the first two estimates is essentially as for
(\ref{th2_3}) and (\ref{th2_4}), and it is given in Section \ref{down}.
As in the system with the wall, we believe these estimates to be of the right order in $L$. 

The estimate (\ref{th5_2}) shows that relaxation in the strictly
delocalized phase is radically different from that of the model with
wall. The proof is based on a somewhat subtle analysis of the behavior
of the signed area under the path. This will be worked out in Section
\ref{52}.  While the logarithmic correction is spurious it might be
that (\ref{th5_2}) captures the correct power law decay of the
spectral gap for $\l<1$, as argued in Section \ref{sec:miclo} below.
Of course, by (\ref{gaptmix}) the bound (\ref{th5_2}) implies that
$\tmix\geq L^{5/2}/(c(\l)(\log L)^8)$.

\subsection{Relaxation of local observables in the localized phase}
\label{sec:relloc}
Finally, we show that in the localized phase local observables
decay to equilibrium following a stretched exponential behavior. For
technical reasons we restrict to the model with the wall.
As it will be apparent from the discussion below, our arguments are
similar to the heuristic ones introduced by D.\ Fisher and D.\ Huse
\cite{Fisher-Huse} in the context of low temperature stochastic Ising
models (see also the more mathematical papers \cite{BoMa} and \cite{FSS}). Specifically,
bounds on the probability of creating an initial local large
fluctuation of the interface around the support of the local function
and on the time necessary in order to make it disappear will play a
key role.

In the localized phase the infinite-volume measure (denoted by
$\mu^+_\infty$) 
is the law of a positive recurrent Markov chain. In order to have more natural statements
in Theorem \ref{th:t13} 
below, we take the
thermodynamic limit as follows. We start from the system with zero
boundary conditions at $\pm L$ for $L\in2\bbN$ (instead of $0,L$
as we did until now) and we denote (with a slight abuse of
notation) by $\mu^{+,\l}_{2L}$ the corresponding equilibrium measure. Then,
for every bounded function $f$ with finite support $\mathcal
S_f\subset \bbZ$, the limit $$
\mu^+_\infty(f):=\lim_{L\to\infty}\mu_{2L}^{+,\l}(f)
$$
exists (cf.\ Lemma \ref{th:tight} and in particular
\eqref{eq:bc}). Similarly, for any fixed $t\geq 0$, 
if $P^+_{t,2L}$ denotes the semigroup 
in the system with zero boundary conditions at $\pm L$, we denote by 
$$
P_t f(\eta):=\lim_{L\to\infty}P^+_{t,2L}f(\eta)\,,
$$
the semigroup
associated to the infinite--volume dynamics in the localized
phase. Standard approximation estimates show that the above pointwise
limit is well defined for every bounded local function $f$
(see e.g.\ the argument in proof of Claim \ref{th:claimpipitilde}
below for more details).  

\begin{Theorem}
\label{th:t13} 
For every $\l>2$ there exists $m>0$ such that the following holds. 


\smallskip

\noindent
1) For every bounded 
local function $f$ there exists a constant $C_f<\infty$ depending on
$\cS_f$
and $\|f\|_\infty$
such that
\begin{eqnarray}
\label{eq:t13}
  \var_{\mu^+_\infty}(P_t f)\leq C_f\, 
\nep{-m\, t^{1/3}},
\end{eqnarray}
for every $t\geq 0$.

\smallskip

\noindent
2) 
For functions $f$ of the form
\begin{eqnarray}
\label{eq:cylinder}
f^{\underline a,I}(\eta):=
1_{\{\eta_x\leq a_x\;\forall x\in I\}} ,
\end{eqnarray}
where $I$ is a finite
subset of $\mathbb Z$ and $a_x\in \mathbb N$, there exists a constant 
$c_f>0$ 
such that
%
%
\begin{eqnarray}
\label{eq:t12}
\var_{\mu^+_\infty}(P_t f)
\geq c_f\,
\nep{- \sqrt{t}/m}\,,
\end{eqnarray}
for every $t\geq 0$.
\end{Theorem}

 The fact that the exponents of $t$ in
(\ref{eq:t13}) and (\ref{eq:t12}) do not match is essentially a
consequence of the fact that the exponents of $L$ in our upper and
lower bounds on the spectral gap in the localized phase also do not match
(cf.\ (\ref{th1_1}) and (\ref{th2_3})).
Theorem \ref{th:t13} 
is proven in Section 
\ref{locco}.

\section{Proof of Theorem \ref{th1} and Theorem \ref{th4}}\la{up}
We are going to use the argument described in Section
\ref{firstarg}. In particular, we recall that both Theorem \ref{th1}
and Theorem \ref{th4} will follow once we show that 
\be
\bbE\widetilde \Phi_t \leq \nep{-\k_L \,t} \widetilde\Phi_0\,,\quad
\;t>0\,,
\la{claims}
\end{equation}
where $\wt \Phi_t$ is given by (\ref{phit}).
Indeed, assuming (\ref{claims}) we can repeat the estimates leading to 
(\ref{arg1}) and (\ref{arg2}) without modifications, which achieves
the proof. 

\subsection{Proof of (\ref{claims}) with the wall}
We shall prove that (\ref{claims}) holds for the system with the wall,
for any $\l>0$. Observe that
\be
\frac{d}{dt}\,\bbE\widetilde \Phi_t = \frac{d}{dt}\,P_t \Phi(\wedge) -
\frac{d}{dt}\,P_t \Phi(\vee) = P_t \cL \Phi(\wedge) - P_t \cL
\Phi(\vee)\,,
\la{cl1}
\end{equation}
where, for simplicity, we omit the $+$ superscript and write
$\cL$ for $\cL^+$ and $P_t$ for $P^+_t$. From Lemma \ref{lemma1} 
and (\ref{arg30}) we
know that 
\be
\cL\Phi = \sum_{x=1}^{L-1}g_x\cL\eta_x = -\k_L\Phi + \Psi\,,
\la{cl2}
\end{equation}
where we use the notation 
\be
\Psi(\eta):=\sum_{x=1}^{L-1}g_x\left[
1_{\{\eta_{x-1} = \eta_{x+1} = 0\}}
- (1-\d) \,1_{\{\eta_{x-1} = \eta_{x+1} =  1\}}\right]\,,
\la{cl3}
\end{equation}
with $\d=2/(1+\l)$. 
Setting $$
\widetilde \Psi_t:=\Psi(\eta^\wedge(t))-\Psi(\eta^\vee(t))\,,$$
equation (\ref{cl1}) becomes
\be
\frac{d}{dt}\,\bbE\widetilde \Phi_t = -\k_L\,\bbE\widetilde \Phi_t
+ \bbE\widetilde \Psi_t\,.
\la{cl4}
\end{equation}
Therefore the claim (\ref{claims}) follows if we can prove that 
\be
\bbE\widetilde \Psi_t \leq 0\,.
\la{cl5}
\end{equation}
It will be convenient to rewrite $ \bbE\widetilde \Psi_t$ as follows.
Define 
\begin{gather*}
  \g_0(x,t) = \bbP(\h^\vee_{x\pm 1}(t)=0) - \bbP(\h^\wedge_{x\pm 1}(t)=0)\,,\\
  \g_1(x,t) = \bbP(\h^\vee_{x\pm 1}(t)=1) - \bbP(\h^\wedge_{x\pm 1}(t)=1)\,.
\end{gather*}
In this way,
\be
\bbE\widetilde \Psi_t = -\sum_{x=1}^{L-1}g_x
\left[\g_0(x,t)-(1-\d)\g_1(x,t)\right]\,.
\la{cl50}
\end{equation}
Clearly, by construction, $\g_0(x,t)=0$ for
$x$ even and $\g_1(x,t)=0$ for $x$ odd.
Note that $\g_i(x,t)\geq 0$ for all $t\geq0$, all $x$ and $i=0,1$,
by monotonicity
(for instance, due to the constraint $\eta_x\ge0$ and to monotonicity 
of the global coupling, $\eta^\vee_{x\pm1}(t)=1$ whenever 
$\eta^\wedge_{x\pm1}(t)=1$, and the non-negativity of $\g_1(x,t)$ immediately
follows). In particular, this implies 
the estimate (\ref{cl5}) if $\l\leq 1$, since in this case $\d\geq 1$. 
The case $\l > 1$ requires more work.

Define $a_x$ as the equilibrium probability that $\h_{x-1}=\h_{x+1}=0$
conditioned to the event that $\h_{x}=\h_{x+2}=1$; similarly, define
$b_x$ as the equilibrium probability that $\h_{x-1}=\h_{x+1}=0$
conditioned on the event that $\h_{x-2}=\h_{x}=1$:
\be
a_x=\mu^+[\h_{x\pm 1}=0\tc \h_{x}=\h_{x+2}=1]\,,\quad 
b_x=\mu^+[\h_{x\pm 1}=0\tc \h_{x-2}=\h_{x}=1]\,.
\la{cl6}
\end{equation}
The proof of (\ref{cl5}) in the case $\l>1$ is based on the next two results.

\begin{Lemma}\label{claimgamma}
For all $t\geq 0$, all $x=2,\dots,L-2$:
\begin{gather}
\g_0(x-1,t)\geq a_{x-1}\g_1(x,t)\,, 
\label{bounds1}\\
\g_0(x+1,t)\geq b_{x+1}\g_1(x,t)\,.
\label{bounds2}
\end{gather}
\end{Lemma}

\begin{Lemma}
\label{th:lemmarho} Set $$\r(x):=\min\{a_{x-1},b_{x+1}\}\,.$$
Then, uniformly in $L$ and $x=2,\ldots,L-2$:
\begin{equation}
\r(x)\geq 1-\d\,.
\label{rx}
\end{equation}
\end{Lemma}

Once we have (\ref{bounds1}) and (\ref{bounds2}) we can estimate
\begin{equation}
\sum_{x=2}^{L-2}g_x\g_1(x,t)\leq \frac12\sum_{x=2}^{L-2}g_x\,
\left\{a_{x-1}^{-1}\,\g_0(x-1,t)+b_{x+1}^{-1}\,\g_0(x+1,t)\right\} \,.
\la{cl7}
\end{equation}
Inserting in (\ref{cl50}) and 
using (\ref{rx}) we arrive at 
\begin{equation}
- \bbE\widetilde \Psi(t)\geq 
\sum_{x=1}^{L-1}\left[g_x-\frac{
g_{x-1} + g_{x+1}}2\right]
\g_0(x,t)\,. 
\label{psi1}
\end{equation}
Recalling 
that $\D g=-\k_L\,g$, 
the
desired claim follows:
$$
- \bbE\widetilde \Psi(t)\geq \k_L
\sum_{x=1}^{L-1}g_x\g_0(x,t) \geq 0\,.
$$

\subsubsection{Proof of Lemma \ref{claimgamma}}
We first prove that for any odd $x=1,\dots,L-3$
\begin{equation}
\bbP(\eta^\wedge_{x-1}(t)=\eta^\wedge_{x+1}(t) = 0)\leq a_{x}\,
\bbP(\eta^\wedge_{x}(t)=\eta^\wedge_{x+2}(t) = 1)\,.
\label{b1}
\end{equation}
Let $A\subset \O_L^+$ denote the subset of non--negative paths $\eta$
such that $\h_{x-1}=\h_{x+1}=0$. Also, let $B\subset \O_L^+$ denote
the subset of non--negative paths $\eta$ such that
$\eta_x=\eta_{x+2}=1$. Note that $A\subset B$. 
If $\mu^+$ denotes the equilibrium measure, we consider the
conditional laws
$\mu_A = \mu^+[\cdot\tc \eta\in A]$ and $\mu_B = \mu^+[\cdot\tc \eta\in
B]$. It is not hard to show that we can find a coupling 
$\nu$ of $(\mu_A,\mu_B)$ such that
$\nu(\h_A\leq \h_B) = 1$ if $\h_A$ is distributed according to
$\mu_A$ and  $\h_B$ is distributed according to
$\mu_B$. As discussed in Section \ref{coupling} this can be obtained 
from the global coupling by letting time go to infinity. 
For any $\xi_A\in A$ and $\xi_B\in B$, we write 
$\nu(\xi_B\tc \xi_A)$
for the $\nu$--conditional probability of having
$\eta_B=\xi_B$ given that $\eta_A=\xi_A$. We have $\nu(\xi_B\tc
\xi_A) = 0$ unless $\xi_B\geq \xi_A$. 

Using the reversibility (\ref{reve}), the left hand side in (\ref{b1}) can
be written as
$$
\sum_{\xi_A\in A}P_t(\wedge,\xi_A) = 
\sum_{\xi_A\in A}P_t(\xi_A,\wedge)\,\frac{\mu^+(\xi_A)}{\mu^+(\wedge)}\,.
$$
Note that for any $\xi_A\leq \xi_B$ monotonicity implies that 
$P_t(\xi_A,\wedge)\leq P_t(\xi_B,\wedge)$. Therefore
we find
\begin{align*}
 \sum_{\xi_A\in A}P_t(\xi_A,\wedge)\,\frac{\mu^+(\xi_A)}{\mu^+(\wedge)} &=
 \sum_{\xi_A\in A} \sum_{\xi_B\in B} \nu(\xi_B\tc \xi_A)
P_t(\xi_A,\wedge)\,\frac{\mu^+(\xi_A)}{\mu^+(\wedge)}\\
&\leq \sum_{\xi_A\in A} \sum_{\xi_B\in B} \nu(\xi_B\tc \xi_A)
P_t(\xi_B,\wedge)\,\frac{\mu^+(\xi_A)}{\mu^+(\wedge)} \\
&= \sum_{\xi_A\in A} \sum_{\xi_B\in B} 
\frac{\nu(\xi_B, \xi_A)}{\mu_A(\xi_A)}
P_t(\wedge,\xi_B)\,\frac{\mu^+(\xi_A)}{\mu^+(\xi_B)}
\end{align*}
Clearly,
$$
\frac{\mu^+(\xi_A)}{\mu_A(\xi_A)} = \mu^+(A)\,,
$$
and
$$
\sum_{\xi_A\in A} 
\nu(\xi_B, \xi_A) = 
\mu_B(\xi_B) = 
\frac{\mu^+(\xi_B)}{\mu^+(B)}\,.
$$
Therefore
$$
\sum_{\xi_A\in A} 
\frac{\nu(\xi_B, \xi_A)}{\mu_A(\xi_A)}
\,\frac{\mu^+(\xi_A)}{\mu^+(\xi_B)} = 
\frac{\mu^+(A)}{\mu^+(B)}
= 
a_{x} \,.
$$
This implies (\ref{b1}).

In a similar way one shows that for any odd $x=1,\dots,L-3$ 
\begin{equation}
\bbP(\eta^\vee_{x-1}(t)=\eta^\vee_{x+1}(t) = 0)\geq a_{x}\,
\bbP(\eta^\vee_{x}(t)=\eta^\vee_{x+2}(t) = 1)\,.
\label{b100}
\end{equation}
The bounds (\ref{b1}) and (\ref{b100}) imply 
(\ref{bounds1}). The complementary
bound (\ref{bounds2}) follows from the same arguments. \qed

\subsubsection{Proof of Lemma \ref{th:lemmarho}} 
We observe that, for $x$ even,
$a_{x-1} = (1-\d/2)\,p_{x}$ where $1-\d/2 = \l/(1+\l)$ is the
equilibrium
probability that
$\h_{x}=0$ given that $\h_{x-1}=\h_{x+1}=1$ and
$p_{x}:=\mu^{+,\l}_{x}(\h_{2}=0)$ is the equilibrium probability that
$\h_{2}=0$ in the system of length $x$.  Similarly, 
$b_{x+1}=(1-\d/2)\,p_{L-x}$. In particular: $$
\r(x)\geq (1-\d/2)\min_{x \;{\rm even}} p_{x}\,.
$$ 
Therefore we need a bound of the form 
\begin{equation}
\min_{x \;{\rm even}}p_{x} \geq \frac{1-\d}{1-\d/2}\,.
\label{bounds3}
\end{equation}
Note that $\frac{1-\d}{1-\d/2} = \frac{\l-1}\l$. We will 
show first that $p_x$ is non-increasing in $x$ and then that
$p_\infty:=\lim_{x\to\infty}p_x\geq (\lambda-1)/\lambda$.  

Indeed, if $x<y$ and $x,y$ are both even, since the function
$f=-1_{\eta_2=0}$ is increasing, the inequality $p_x\geq p_y$ is
easily derived from the FKG inequality (\ref{mono1}). Next, note that
from \eqref{eq:id2} $$
p_\infty=\lim_{x\to\infty}\mu^{{\l}/2}_x(\eta_2=0) $$ where $\mu^\l_x$
denotes the equilibrium measure without the wall in the system of
length $x$ and with parameter $\l$.  From \cite[Th. 2.3]{cf:GB} (or
from Theorem \ref{th:asymptZ} above) we have
\begin{eqnarray}
\lim_{x\to\infty}\mu^{\l}_x(\eta_2=0)=
\lambda/2\exp(-2\, F(\lambda))
\end{eqnarray}
where $F$ is the free energy defined in Section \ref{sec:eq}.  Since
$F(\l)=0$ for $\l\leq 1$ we have $p_\infty=\lambda/4\geq (\l-1)/\l$
for all $\l\leq 2$.  As for $\l\geq 2$, one uses the explicit
expression
\eqref{eq:tildeFl1}, which gives $p_\infty=(\l-1)/\l$.
This ends the proof. \qed

\subsection{Proof of (\ref{claims})  without the wall}
Here we assume $\l\geq 1$. In the model without the wall, 
we can repeat the computations leading to
(\ref{cl4}). The function $\Psi$ containing the correction terms
from Lemma \ref{lemma1} is now given by 
\be
\Psi(\eta):=(1-\d)\sum_{x=1}^{L-1}g_x\left[
1_{\{\eta_{x-1} = \eta_{x+1} =- 1\}}
- \,1_{\{\eta_{x-1} = \eta_{x+1} =  1\}}\right]\,,
\la{cl33}
\end{equation}
with $\d=2/(1+\l)\leq 1$. 
Setting again $\widetilde
\Psi_t:=\Psi(\eta^\wedge(t))-\Psi(\eta^\vee(t))\,,$
we arrive at the same expression given in (\ref{cl4}). 
Therefore it suffices to show 
that $\bbE\wt\Psi_t \leq 0$. This in turn is an immediate consequence
of the next lemma. 

\begin{Lemma}\la{nowall1}
For every even $x$ we have 
\begin{gather}
\bbP(\h^\wedge_{x\pm 1}(t)=1)\geq \bbP(\h^\wedge_{x\pm
  1}(t)=-1)\,,\la{cl44}\\
\bbP(\h^\vee_{x\pm 1}(t)=-1)\geq \bbP(\h^\vee_{x\pm
  1}(t)=1)\,.\la{cl45}
\end{gather} 
\end{Lemma}
\proof
By symmetry it suffices to prove (\ref{cl44}) only.
We use an argument similar to that of Lemma \ref{claimgamma}.
Namely, call $A$ the set of all paths $\eta\in\O_L$ such that
$\h_{x\pm 1}=-1$ and 
 $B$ the set of all paths $\eta\in\O_L$ such that
$\h_{x\pm 1}=1$ and 
let
$\mu_A=\mu[\cdot\tc\eta\in A]$, $\mu_B=\mu[\cdot\tc\eta\in B]$.
Again, using the global coupling 
we construct a coupling 
$\nu$ of $(\mu_A,\mu_B)$ such that
$\nu(\h_A\leq \h_B) = 1$ if $\h_A$ is distributed according to
$\mu_A$ and  $\h_B$ according to
$\mu_B$. 
We have $\nu(\xi_B\tc
\xi_A) = 0$ unless $\xi_B\geq \xi_A$. Using monotonicity we also have 
$P_t(\xi_A,\wedge)\leq P_t(\xi_B,\wedge)$ whenever $\xi_B\geq \xi_A$.
Therefore the same computation as in the proof of Lemma
\ref{claimgamma} now gives
\begin{align*}
\bbP(\h^\wedge_{x\pm
  1}(t)=-1) &= \sum_{\xi_A\in A} P_t(\wedge,\xi_A) 
\leq \frac{\mu(A)}{\mu(B)} \sum_{\xi_B\in B} P_t(\wedge,\xi_B) =  
\bbP(\h^\wedge_{x\pm 1}(t)=1)\,,
\end{align*}
where we use the symmetry $\mu(A)=\mu(B)$. \qed

\section{Proof of Theorem \ref{th2} and related bounds}\la{down}

\subsection{Upper bounds on the spectral gap}
\label{sec:gap1L}
We start with the proof of (\ref{th2_3}). Note that this bound can be derived
from the independent estimate (\ref{th2_4}) by using (\ref{gaptmix}). However,
we show an explicit test function which
reproduces the bound $\gap \leq c/L$ in the localized phase (\ie when $\l>2$
for $\cL^+$ and when $\l>1$ for $\cL$). The idea is reminiscent of 
an argument used in \cite{BoMa} for the low temperature Ising model. 
\begin{Proposition}
\label{th:gapL1}
For every $\l>2$ there exists a constant $c\in(0,\infty)$
such that
\begin{eqnarray}
\frac{ {\mathcal E}^+(f_c,f_c)}
{\var(f_c)}\leq \frac{32c^2}{L}
\end{eqnarray}
where for $a>0$
\begin{eqnarray}
f_a(\h):=\exp\left(\frac {a}L \sum_{x=1}^{L-1} \h_x\right).
\end{eqnarray}
The same holds without the wall, if $\l>1$.
\end{Proposition}

{\sl Proof of Proposition \ref{th:gapL1}}.
Let
\begin{eqnarray}
\psi(a):=\lim_{L\to\infty}\frac1L \log \mu^{+,\l}_L(f_a).
\end{eqnarray}
The limit exists since 
$$
\log \mu^{+,\l}_L(f_a)=\log
\sum_{\eta\in \Omega^+_L} \left(\lambda^{N(\eta)}f_a(\eta) \right)-
\log Z^+_L(\l)
$$
and both terms in the right-hand side are super-additive in $L$
(both terms, once divided by $L$, tend to a finite limit for
$L\to\infty$: for the second one this follows from Theorem \ref{th:asymptZ}, 
while for the first one just note that $\lambda^{N(\eta)}f_a(\eta) 
\le\l^L e^{aL}$). 
Observe also that
$\psi(a)$ is non-decreasing in $a$: this is obvious in presence of 
the wall since $\eta_x\geq0$, and in absence of the wall
this follows from the fact that $\psi(a)$ is convex and that
$$
\left.\partial_a \log \mu^\l_L(f_a)\right|_{a=0}=0.
$$
 We will show that
there exists $\bar a\in(0,\infty)$ such that $\psi(a)=0$ for $a\leq
 \bar a$ and $\psi(a)>0$ for $a> \bar a$. Then, choosing $
\bar a/2<c<\bar a$, one has that, for $L$ sufficiently large,
\begin{eqnarray}
\var(f_c)\geq (1/2)\mu^{+,\l}_L(f_{c}^2).
\end{eqnarray}
As for the Dirichlet form, since 
\begin{eqnarray}
\left |Q^+_x f_c(\h)-f_{c}(\eta)\right|\leq \frac{4c}{L} 
f_{c}(\eta),
\end{eqnarray}
one deduces easily from \eqref{eq:dirich} that 
\begin{eqnarray}
{\mathcal E}^+(f_{c},f_{c})\leq \frac{16c^2}L 
\mu^{+,\l}_L(f_{c}^2).
\end{eqnarray}
The statement of the proposition is then proven once we show the
existence of the $\bar a$ introduced above, and it is here that the
assumption $\l>2$ will play a role. Note that $\bar a$ is uniquely
defined by the monotonicity of $\psi(a)$, so we have only to show that
it is neither zero nor infinity. First of all,
\begin{eqnarray}
\mu^{+,\l}_L\left(f_a
\right)\geq \mu^{+,\l}_L\left(f_a\, 1_{\eta=\wedge}\right)\geq e^{(a/8)L}
\end{eqnarray}
if $a$ is sufficiently large, so that $\bar a<\infty$. Conversely, 
from Jensen's inequality $$f_a\leq \frac1L\sum_{i=1}^L \nep{a\eta_x}$$
and from Lemma \ref{th:tight} we see that $\mu^{+,\l}_L(f_a)\leq C$ for some 
$C$ independent of $L$ if, say, $a<F(\l)/2$. In this case $\psi(a)=0$
and therefore $\bar a>0$.

The proof of the Proposition in the case of $\l>1$ and no wall is
essentially identical. \qed

\smallskip

In the delocalized phase and at the critical point,
on the other hand, we can do better and show a test function 
which gives the (optimal) behavior of order $L^{-2}$. The next result proves
Eq. \eqref{th2_1}. 

\begin{Proposition}
\label{th:gapL2}
Let 
\begin{eqnarray}
f(\eta):=\sum_{x=1}^{L}
\sin\left(\frac{\pi(x-1/2)}L\right)(\eta_x-\eta_{x-1}).
\end{eqnarray}
There exists $c<\infty$  such that for every $\l\leq2$ 
 the following holds:
\begin{eqnarray}
\frac{ {\mathcal E}^+(f,f)}{\var(f)}\leq \frac{c}{L^2}.
\end{eqnarray}
\end{Proposition}

\proof
Set for convenience
$\phi_x:=\h_x-\h_{x-1}$ and $h_L(x):=\sin(\pi(x-1/2)/L)$. It is clear 
from the symmetry $x\leftrightarrow (L-x)$ 
that $\mu_{L}^{+,\l}(f)=0$.  The Dirichlet form is easily
estimated from above: since the occurrence of a flip at site $x$ has
the effect of exchanging the values of $\phi_x$ and $\phi_{x+1}$, one
has
\begin{eqnarray}
{\mathcal E}^+(f,f)\leq 
\sum_{x=1}^L 
\mu^{+,\l}_L(\phi_x\ne\phi_{x-1})
\left[h_L(x+1)-h_L(x)\right]^2.
\end{eqnarray}
 Therefore, for 
$L$ sufficiently large,
one has
\begin{eqnarray}
\label{eq:UBdirichelet}
{\mathcal E}^+(f,f)\leq \frac{c}L \int_0^1 \cos\left(\pi s\right)^2
\,ds\,.
\end{eqnarray}
As for the second moment of $f$, 
what we need to show is that 
\begin{eqnarray}
\label{eq:itsuffices}
 \mu_L^{+,\l}(f^2) \geq c L
\end{eqnarray}
with $c\in(0,\infty)$.

To this end we note that
\begin{eqnarray}
L^{-1/2}f(\h)&=&\frac1{\sqrt L}\sum_{x=1}^{L-1}\h_x\left(
h_L(x)-h_L(x+1)
\right)\\\nonumber
&=&-\frac\pi{ L^{3/2}}\sum_{x=1}^{L-1}\h_x
\cos\left(\frac{\pi(x-1/2)}L\right)+O(L^{-1/2})
\end{eqnarray}
where the estimate on the error term is uniform in $\h$.
Introducing the continuous, piecewise linear process
$\{\eta^{(L)}_s\}_{s\in[0,1]}$
such that $\eta^{(L)}_{x/L}:=L^{-1/2}\h_x$ for $x=0,\ldots,L$ 
and $\partial^2_s\eta^{(L)}_s=0$
for $s\in (x/L, (x+1)/L)$,
one has
\begin{eqnarray}
&&L^{-1/2}f(\h)
=- \pi\int_0^1ds\, \eta^{(L)}_s \cos(\pi s)+O(L^{-1/2}).
\end{eqnarray}
It is easy to see that for every $k>0$ 
\begin{eqnarray}
\label{eq:itisisi}
\sup_L \sup_{x\in \L} \mu_L^{+,\l}\left(
\frac{(\h_x)^k}{L^{k/2}}\right)
<\infty.
\end{eqnarray}
(Indeed, using monotonicity a couple of times,
$$
\mu_{L}^{+,\l}\left(\frac{(\h_x)^k}{L^{k/2}}\right )
\leq
\bE\left(\left.\frac{(\h_{L})^k}{L^{k/2}}\right|\h_y>0\,,\,\; \forall 1\leq 
y\leq L\right)
$$
and the latter expression is 
seen to be finite uniformly in $L$ 
using the ``Ballot Theorem'' \cite[Sec. III.1]{cf:Feller} plus 
Stirling's formula.)
Equation \eqref{eq:itsuffices} is therefore proven if we show that
\begin{eqnarray}
\label{eq:almost}
\int_0^1 dt\int_0^1ds\,\cos(\pi t)\cos (\pi s)
\mu_L^{+,\l}\left(\h^{(L)}_t\h^{(L)}_s\right)
\stackrel{L\to\infty}{\to}
c\in(0,\infty).
\end{eqnarray} 
Consider first the case $\l<2$, in which case
\begin{eqnarray}
\label{eq:convcovar}
\mu_L^{+,\l}\left(\h^{(L)}_t\h^{(L)}_s
\right)\stackrel{L\to\infty}{\to} \mu_X(X_t X_s),
\end{eqnarray}
where $X_.$ is the Brownian bridge in $[0,1]$ conditioned to be
non-negative (\ie, the Bessel bridge of dimension $3$ from $0$ to
$0$), whose law we denote by $\mu_X(\cdot)$. Equation
\eqref{eq:convcovar} follows from the convergence in distribution  of 
the process $\h^{(L)}_.$ to 
$X_.$ (this is proven in \cite{cf:DGZ} in a slightly different
setting, but the techniques developed there can be extended
to our case; see also \cite{Iso_Yos}), together with the uniform
integrability \eqref{eq:itisisi}.  Thanks to \eqref{eq:itisisi} we can
take the $L\to\infty$ limit inside the integral in
\eqref{eq:almost} and the limit is
$$
\var_{\mu_X}\left(\int_0^1 ds\, \cos(\pi s)X_s\right)>0.
$$
 Together with \eqref{eq:UBdirichelet}
this concludes the proof of the proposition in the strictly delocalized
case. 

It remains to consider the critical case $\l=2$.
In this case, \cite[Th. 5.1]{cf:CGZ} plus \eqref{eq:itisisi} imply that
\begin{eqnarray}
\mu_L^{+,2}\left(\h^{(L)}_t\h^{(L)}_s
\right)\stackrel{L\to\infty}{\to} \mu_B(|B_t|\, |B_s|),
\end{eqnarray}
$B_.$ being the Brownian bridge on $[0,1]$ with law $\mu_B$.
Therefore, one finds in this
case 
\begin{eqnarray}
\label{eq:weakconv}
\mu^{+,2}_L(f^2/L)\stackrel{L\to\infty}\longrightarrow\pi^2
\var_{\mu_B}\left(\int_0^1ds\,\cos(\pi s)|B_s|\right)>0
\end{eqnarray}
and the conclusion follows as before. \qed 

\subsection{Lower bounds on mixing times}
We begin with the proof of (\ref{th2_2}) for the 
system with the wall at $0<\l <2$.

\begin{Proposition}
\la{loboth1}
For any $\l<2$, 
\be
\tmix\geq \left(\frac1{2\pi^2} + o(1)\right)\,L^2\log L\,.
\la{lobo1}
\end{equation}
\end{Proposition}
Let $\Phi$ and $\Psi$ be the functions defined in (\ref{phi}) and
(\ref{cl3}) respectively.  We consider their evolutions when
the process starts in the maximal configuration $\wedge$ and 
set $\hat\Phi_t: =
\Phi(\eta^\wedge(t))$ and $\hat\Psi_t: =
\Psi(\eta^\wedge(t))$. 
From the computation in Lemma
\ref{lemma1}, as in (\ref{cl2}) we obtain
\be
\frac{d}{dt}\bbE[\hat \Phi_t] = P_t\cL\Phi(\wedge)=
-\k_L\,\bbE[\hat \Phi_t] + \bbE[\hat \Psi_t]\,,
\la{lobo3}
\end{equation}
where as usual we omit the $+$ superscript and write
$\cL$ for $\cL^+$ and $P_t$ for $P^+_t$.
Next, we claim that for any $\l<2$, there exists $c(\l)<\infty$ 
such that for all $t\geq 0$, $L\geq 2$: 
\be
\bbE[\hat \Psi_t] \geq - c(\l)\, L^{-1/2}\,.
\la{lobo4}
\end{equation}
To prove this, observe that 
$$
\bbE[\hat\Psi_t]\geq -(1-\d)\sum_{x=1}^{L-1}g_x
\bbP(\eta^\wedge_{x\pm 1}(t)=1)\,.
$$
Note that if $\l\leq 1$ then $\d\geq 1$ and therefore $\bbE[\hat\Psi_t]\geq 0$.
If $\l\in (1,2)$ we use the following argument to prove (\ref{lobo4}).
Monotonicity allows us to bound 
$\bbP(\eta^\wedge_{x\pm 1}(t)=1)$ from
above by the equilibrium probability
$\mu^+(\eta_{x\pm 1}=1)$. The latter, in turn, for each
$2\leq x\leq L-2$ and $\l<2$ 
is estimated with 
\be
\mu^+(\eta_{x\pm 1}=1) = \frac{1+\l}{\l}\,\mu^+(\eta_x=0)
\leq c(\l)\,\frac{L^{3/2}}{(L-x)^{3/2}x^{3/2}}\,,
\la{lobo5}
\end{equation}
where $c(\l)$ is a suitable constant.
Note that (\ref{lobo5}) follows from (\ref{eq:asymptZ}) using 
$$
\mu^+(\eta_x=0) = \frac{Z^+_x(\l)Z^+_{L-x}(\l)}{Z^+_L(\l)}
=\frac12\, \frac{Z_x(\l/2)Z_{L-x}(\l/2)}{Z_L(\l/2)}\,,
$$
and the fact that $\l<2$. 
Once we have (\ref{lobo5}), using $g_x\leq \frac{\pi x}L$ for $x\leq
L/2$ and $g_x\leq \frac{\pi (L-x)}L$ for $x\geq L/2$ we obtain
$$
\bbE[\hat\Psi_t]\geq
-c(\l)\,\sum_{x=1}^{L/2}\frac{1}{x^{1/2}\,L}\,,
$$
with a new constant $c(\l)$. This implies the claim (\ref{lobo4}).

\smallskip
Next, we integrate (\ref{lobo3}) using (\ref{lobo4}) to obtain
\begin{align*}
\bbE[\hat \Phi_t] &\geq \nep{-\k_L t}\hat \Phi_0 -
c(\l)\,L^{-1/2}\int_0^t\nep{-\k_L(t-s)}\,ds\, \\
&\geq 
\nep{-\k_L\,t}\hat \Phi_0 - c(\l) \,L^{-1/2}\,\k_L^{-1}\,.
\end{align*}
Therefore we have shown that for each $\l<2$, for some constant $c(\l)$, 
for all $t\geq 0$ and $L\geq 2$:
\be
\bbE[\hat \Phi_t]\geq \nep{-\k_L t}\hat \Phi_0 - c(\l)\,L^{3/2}\,.
\la{lobo6}
\end{equation}
Since $\hat \Phi_0\geq c\, L^2$ and $\k_L\sim\pi^2/2L^2$, 
from (\ref{lobo6}) we see that
$\bbE[\hat \Phi_t]$ is much larger than its
equilibrium value $\bbE[\hat \Phi_\infty]=O( L^{3/2})$ for times $t$
within, say, $\frac{1}{2\pi^2}\,L^2\log L$. 
However, this is still not enough to prove that the mixing
time is at least of order $L^2 \log L$, since the $L_\infty$ norm of
$\Phi$ is of order $L^2$. 

\smallskip

Following Wilson
\cite{cf:Wilson}, we turn to an estimate on the variance of $\hat \Phi_t$.
\begin{Lemma}
\la{levar}
For every $\l<2$ there exists $c(\l)<\infty$ such that for all $t>0$:
\be
\var(\hat \Phi_t) = 
\bbE[\hat \Phi_t^2] - \bbE[\hat \Phi_t]^2 \leq c(\l) \,L^{7/2}\,.
\la{lobo12}
\end{equation}
\end{Lemma}
\proof 
We start by giving an upper bound on $\bbE[\hat \Phi_t^2]$.
Recall from (\ref{gene}) that 
$$
\cL\Phi^2 = \sum_{x=1}^{L-1} 
\left [Q_x\Phi^2 - \Phi^2\right]\,, 
$$
where $Q_x$ is the equilibrium measure at $x$ conditioned on the
configuration $\eta$ outside of $x$. Writing 
$Q_x(\xi\tc \eta):=\mu^+[\xi\tc \eta_y\,,\;y\neq x]$ for the
associated kernel, for every
$\eta\in\O^+_L$ we have
\begin{align*}
Q_x\Phi^2(\eta) &= \sum_{\xi} 
Q_x(\xi\tc \eta)
\,\Phi^2(\xi)
\\&= \sum_{\xi} 
Q_x(\xi\tc \eta)
\left[\Phi^2(\eta) +
  2\Phi(\eta)(\Phi(\xi)-\Phi(\eta))+ (\Phi(\xi)-\Phi(\eta))^2\right]\,.
\end{align*}
We can estimate $\sum_{\xi} 
Q_x(\xi\tc \eta)(\Phi(\xi)-\Phi(\eta))^2\leq 4$
for each $x$ and $\eta$. Indeed, each transition can at most change
the function $\Phi$ by $2$.  
Therefore
\be
\cL\Phi^2(\eta)\leq 4L + 2\Phi(\eta)\cL\Phi(\eta)
= 4L - 2\k_L\Phi(\eta)^2 + 2\Phi(\eta)\Psi(\eta)\,,
\la{lobo7}
\end{equation}
where we have used again (\ref{cl2}) in the last step. 
In conclusion, inserting (\ref{lobo7}) in the identity
$$
\frac{d}{dt}\bbE[\hat \Phi_t^2] = P_t
\cL\Phi^2(\wedge)\,,$$ 
and integrating we obtain 
\be
\bbE[\hat \Phi_t^2] \leq 
\nep{-2\k_L t}(\hat \Phi_0)^2 + 4 L
\k_L^{-1} + 2\int_0^t\bbE[\hat \Phi_s\hat \Psi_s]\nep{-\k_L(t-s)}ds\,.
\la{lobo8}
\end{equation}
To estimate the last term in (\ref{lobo8}) we note that
$\hat \Phi_s\leq \hat\Phi_0 = O(L^2)$ uniformly. Moreover, for any
$\l>0$, $s\geq 0$: 
\begin{align*}
\bbE[\hat\Psi_s]&\leq 
\sum_{x=1}^{L-1}g_x\left[
\bbP(\eta^\wedge_{x\pm 1}(s)=0)
+ \bbP(\eta^\wedge_{x\pm 1}(s)=1)\right]\\
&\leq \sum_{x=1}^{L-1}g_x\left[\mu^+(\eta_{x\pm 1}=0)+
\mu^+(\eta_{x\pm 1}=1) \right]\,,
\end{align*}
where the last step follows from monotonicity. 
As in (\ref{lobo5}) we have the following 
equilibrium bounds valid for any $\l<2$:
\be
\mu^+(\eta_{x\pm 1}=i)
\leq c(\l)\,\frac{L^{3/2}}{(L-x)^{3/2}x^{3/2}}\,,\quad i=0,1\,.
\la{lobo9}
\end{equation}
As in the proof of (\ref{lobo4}) these estimates imply
\be
\bbE[\hat\Psi_s]\leq c(\l)\,L^{-1/2}\,.
\la{lobo10}
\end{equation}
From (\ref{lobo8}) we therefore obtain
\begin{align}
\bbE[\hat\Phi_t^2] &\leq \nep{-2\k_L t}\hat \Phi_0^2 + 4 L
\k_L^{-1} + c(\l) \k_L^{-1}\hat \Phi_0 \,L^{-1/2}\nonumber\\
& \leq \nep{-2\k_L t}\hat \Phi_0^2 +
c(\l)L^{7/2}\,.
\la{lobo11}
\end{align}
From (\ref{lobo6}) we know that 
$\bbE[\hat \Phi_t]^2\geq \nep{-2\k_L t}\hat \Phi_0^2 -
c(\l)\, L^{7/2}$ so that we deduce the upper bound (\ref{lobo12}).
\qed

\smallskip \smallskip \smallskip 
Using Lemma \ref{levar} we can finish the proof of Proposition \ref{loboth1}.
Letting $t\to\infty$ in (\ref{lobo12}) we obtain a bound
on the equilibrium variance
\be
\var_{\mu^+}(\Phi) = \var(\hat \Phi_\infty) \leq 
c(\l) \,L^{7/2}\,.
\la{lobo13}
\end{equation}
Define the set 
$$
A_\g = \{\eta:\;\Phi(\eta)\leq L^{2-\g}
\}\, ,\quad \g\in\left(0,1/4\right)\,. 
$$
Since $\mu^+(\Phi)\leq \sum_{x}\mu^+(\eta_x)\leq c\, L^{3/2}$ we
see that, from Chebyshev's inequality and (\ref{lobo13}):
\begin{align*}
1- \mu^+(A_\g) &\leq 
\mu^+\left(|\Phi-\mu^+(\Phi)|\geq \frac12
  \,L^{2-\g}\right)\\& 
\leq 4\,L^{-4+2\g}\var_{\mu^+}(\Phi) \leq
c(\l)\,L^{-\frac12 + 2\g}\,.
\end{align*}
Let $P_t(\wedge,\cdot)$ denote the distribution of 
$\eta^\wedge(t)$. Using (\ref{lobo6}) we see that if 
$t\k_L \leq a \log L$ for some $a<\g$ then $\hat \Phi_t\leq L^{2-\g}$ implies 
$|\hat \Phi_t-\bbE(\hat \Phi_t)|\geq c \,L^{2-a}$, for some $c>0$, 
for all $L$ large enough. 
From Chebyshev's inequality and (\ref{lobo12}) we then have
\begin{align*}
P_t(\wedge,A_\g) &
\leq \bbP\left(|\hat \Phi_t-\bbE(\hat \Phi_t)|\geq c
  \,L^{2-a}\right)\\& \leq c^{-2}
\,L^{-4+2a}\var(\hat \Phi_t) \leq
c(\l)\,L^{-\frac12 + 2a}\,.
\end{align*}
In conclusion, taking $\g=(\frac14 - \e)$, 
$a=\g-\e$ we see that
for $L$ sufficiently large we have
\be
\|P_t(\wedge,\cdot)-\mu^+\|_{\rm var} \geq 
|P_t(\wedge,A_\g)-\mu^+(A_\g)|\geq  1 - L^{-\e}\,,
\la{lobo130}
\end{equation}
whenever $t\leq \left(\frac14 -2\e\right)\k_L^{-1}\log L$.  Since
$\k_L\sim\pi^2/2L^2$, this ends the proof of Proposition
\ref{loboth1}. \qed

\newcommand{\frT}{\ensuremath{\mathfrak T}}
\subsection{A universal lower bound on the mixing time}
\label{sec:universal}
Here we shall prove the bound (\ref{th2_4}) and the corresponding
estimate in Theorem \ref{th5}.
\begin{Theorem}
\la{unilobo}
Both with and without the wall, for every $\l>0$:
\be
\tmix \geq \,L^2/32\,.
\la{lobo18}
\end{equation}
\end{Theorem}
The proof is divided in three steps. First we prove the
statement at $\l=\infty$ for the system with the wall. Then we extend
it to any $\l>0$ with the wall and finally we show how to prove it
for all $\l>0$ without the wall.

\subsubsection{$\l=\infty$ with the wall} 
Recall that here $\mu^+=\d_{\vee}$ is the Dirac mass at the minimal
configuration $\vee$.  
In this situation, the definitions of generator, spectral gap and
Dirichlet form as given in Section \ref{sec:dyn} do not make sense,
and the dynamics is defined as in the introduction, with $\l/(1+\l)$
replaced by $1$. In other words, the rules for updating $\h_x$
become:
\begin{itemize}
\item if $\eta_{x-1}\neq \eta_{x+1}$, do nothing; 
\item if 
$\eta_{x-1}= \eta_{x+1}=j$ and $j\neq 1$, set $\eta_x = j\pm 1$
with equal probabilities; 
\item if $\eta_{x-1}= \eta_{x+1}= 1$, set $\eta_x =0$ with
probability 1.
\end{itemize}
Similar considerations hold for the $\l=\infty$ dynamics without the wall.

We want to estimate the expected time needed to
go from $\wedge$ to $\vee$.  Following a well known argument using the
mapping with simple exclusion (see e.g.\ the proof of Theorem 1.3 in
\cite{FSS}) we would obtain a bound of the form $L^2/\log L$. We shall
remove the spurious $\log L$ factor by means of the following
argument.

\begin{figure}[h]
\centerline{
\psfrag{a}{$d_1$}
\psfrag{b}{$d_2$}
\psfrag{c}{$d_3$}
\psfrag{d}{$d_4$}\psfrag{q1}{$q_1$}\psfrag{q2}{$q_2$}
\psfig{file=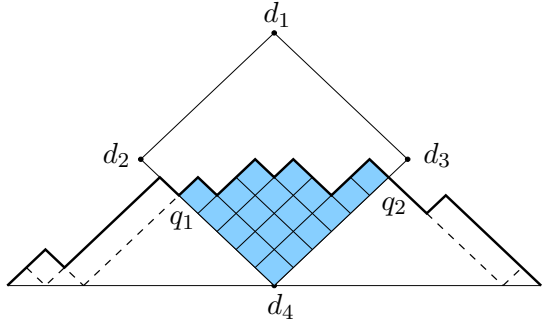,height=1.59in,width=2.8in}
}
\caption{The square $D$ and the area $B_D(\eta)$ for a given path $\eta\in\O_L^+$.} 
\label{fig30}
\end{figure}

We suppose for simplicity that $L/2$ is even (it is straightforward to
modify the construction in the case $L/2$ odd).
Let $D$ denote the square identified by the four vertices with
coordinates
$d_1=(\frac{L}2,\frac{L}2)$, 
$d_2=(\frac{L}4,\frac{L}4)$, 
$d_3=(\frac{3L}4,\frac{L}4)$ and $d_4=(\frac{L}2,0)$.
Given a path $\eta\in\O_L^+$ we call 
$q_1(\eta)$ and $q_2(\eta)$ the points where $\eta$ crosses the lines
$d_2-d_4$ and $d_3-d_4$, respectively, with the rule
that if the path $\eta$ touches the line in more than one point we use
the lowest, \ie the closest to $d_4$, see Figure \ref{fig30}. 
We call $\eta_D$
the portion of the path $\eta$ between  $q_1(\eta)$
and  $q_2(\eta)$. Also, we need to introduce the map which
associates each path $\eta$ with the minimal path compatible with the
portion $\eta_D$, see the dashed lines in Figure \ref{fig30}. 
We write $\frT_D\eta$ for this new configuration.
Note that $\frT_D\eta\leq \eta$ for every $\eta$. If $\eta$ is
such
that $q_1(\eta)=q_2(\eta)=d_4$, \ie if $\eta_{\frac{L}2}=0$, then
$\frT_D\eta = \vee$.
Also, 
$\frT_D\wedge = \wedge$. 

Let $B_D(\eta)$ denote the area of the region
inside the square $D$ enclosed by the path $\eta_D$ and the broken line 
joining points $q_1(\eta),d_4,q_2(\eta)$. 
The area $B_D$ will be measured by the number
of 
elementary $\sqrt 2 \times \sqrt 2$ squares it contains, so that \eg $B_D(\wedge)=L^2/16$. 

\smallskip

\smallskip

Let $\eta(t)=\eta^\wedge(t)$ denote the time evolution starting from
$\eta(0)=\wedge$. We shall consider a modified evolution $\xi(t)$ that
can be coupled to $\eta(t)$ in such a way that $\xi(t)\leq
\eta(t)$. We use the same Poisson clocks and the same ``coins'' (see
Section \ref{coupling}) for the two processes. When a clock rings we
update as in the global coupling with
$(\eta,\xi)\to(\eta',\xi')$. After this updating the configuration
$\xi'$ is replaced by $\frT_D\xi'$. The net result is therefore the
updating $(\eta,\xi)\to(\eta',\frT_D\xi')$ where
$(\eta,\xi)\to(\eta',\xi')$ is a standard update under the global
coupling. Since $\frT_D\xi\leq \xi$ we have $\eta(t)\geq \xi(t)$
almost surely.

We set $\xi(0)=\eta(0)=\wedge$. We want to estimate the expected value
of $ B_D(\xi(t))$ from below. Observe that $\xi(t)$ lives in the space
$\wt\O=\{\si\in
\O_L^+:\, \frT_D\si=\si\}$. If  $\cG$ stands for the generator of the Markov
chain $\xi(t)$, then we claim that, for any $\xi\in\wt\O$ 
\be
\cG B_D(\xi) \geq -1\,.
\la{lobo21}
\end{equation}
To prove (\ref{lobo21}) note that $B_D(\xi)$ can only change by $\pm1$
according to whether there is a mountain/valley in the path
$\xi_D$. 
Each valley in $\xi_D$ contributes with $\frac12$ to
(\ref{lobo21}), if $\xi\neq \vee$. Moreover, each mountain in $\xi_D$
contributes with $-\frac12$ to (\ref{lobo21}) unless reversing it
would result in the configuration $\vee$, in which case its
contribution to (\ref{lobo21}) is $-1$.  However the number of
mountains minus the number of valleys in $\xi_D$ is always $1$ (unless
$\xi=\vee$, in which case $\cG B_D(\xi)=0$). This implies
(\ref{lobo21}).

From (\ref{lobo21}) we know that the martingale 
$$M_t = B_D(\xi(t)) -
B_D(\xi(0)) - \int_0^t\cG B_D(\xi(s)) ds\,,$$ satisfies $0 = \bbE M_t
\leq \bbE B_D(\xi(t)) - B_D(\xi(0)) + t$ or, 
\be
\bbE B_D(\xi(t))\geq  B_D(\wedge) -\, t \,.
\la{lobo22}
\end{equation}
Setting $f(\eta):=B_D(\eta)/B_D(\wedge)$,
we have $\bbE f(\eta(t)) \geq \bbE B_D(\xi(t))/B_D(\wedge)$
and $\d_{\vee}(f)=0$. Moreover,
$\|f\|_\infty\leq 1$ and therefore
\begin{align}
\|P_t(\wedge,\cdot)-\d_{\vee}\|_{\rm var} &
\geq |P_t(\wedge,\cdot)(f)-
\d_{\vee}(f)| \nonumber\\
& = \bbE f(\eta(t)) 
\geq 1 - \frac{t}{B_D(\wedge)}\,.
\la{lobo23}
\end{align}
Since $B_D(\wedge)=\frac{L^2}{16}$ 
we have $\tmix \geq (\nep{}-1)L^2/16\nep{}\geq L^2/32$. \qed

\subsubsection{$\l > 0$ with the wall}
Let $\eta(t)$ denote the evolution with the wall, for a given $\l>0$,
and with maximal initial condition $\eta(0)=\wedge$. If $\xi(t)$ denotes
the process defined above, with $\xi(0)=\wedge$, we can couple the two
processes in such a way that $\eta(t)\geq \xi(t)$ almost surely
and therefore $\bbE B_D(\eta(t))\geq \bbE B_D(\xi(t))\geq B_D(\wedge)
-t$
by (\ref{lobo22}).  
Set again $f(\eta):=B_D(\eta)/B_D(\wedge)$.
Since
$\mu^+(B_D)\leq \mu^+(A) = O(L^{3/2})$ uniformly in $\l>0$, the
equilibrium average of $f$ satisfies 
$\mu^+(f) = O(L^{-1/2})$. Therefore
\begin{align}
\|P_t(\wedge,\cdot)-\mu^+\|_{\rm var} &
\geq |P_t(\wedge,\cdot)(f)-
\mu^+(f)| \nonumber\\
& = \bbE f(\eta(t)) + O(L^{-1/2})
\geq 1 - \frac{t}{B_D(\wedge)} + O(L^{-1/2})\,.
\la{lobo29}
\end{align}
As in (\ref{lobo23}) we obtain 
$\tmix \geq L^2/32$ provided $L$ is sufficiently large. \qed

\subsubsection{$\l > 0$ without the wall}
Call $\eta(t)$ the evolution without the wall, for a given $\l>0$, and
with maximal initial condition $\eta(0)=\wedge$. We can use the same
arguments given above but we have to modify the process $\xi$ in order
to satisfy the monotonicity $\eta(t)\geq \xi(t)$. Recall the
construction of the square $D$ and the associated path $\eta_D$, see
Figure \ref{fig30}. The transformation $\frT_D$ here will be defined
as follows. Given the portion of the path $\eta_D$ then $\frT_D\eta$
is the minimal configuration $\eta'\in\O_L$ (\ie without the wall)
such that $\eta'_D=\eta_D$. Also, we add the rule that if
$\eta_{\frac{L}2}\leq 0$ then $\frT_D\eta=\vee = -\wedge$.  In this
way, if the process $\xi(t)$ is defined as before (but with the new
$\frT_D$), then we can guarantee the domination $\eta(t)\geq
\xi(t)$. In particular, $B_D(\eta(t))\geq B_D(\xi(t))$.  Note that, by
the same arguments, our process $\xi(t)$ satisfies (\ref{lobo21}) and
therefore (\ref{lobo22}). Then we set $f(\eta):=B_D(\eta)/B_D(\wedge)$
and observe that the equilibrium average of $f$ satisfies $\mu(f) \leq
\mu^+(f) = O(L^{-1/2})$ uniformly in $\l>0$. The rest of the argument
is the same as for (\ref{lobo29}). \qed

\subsection{On the mixing time at $\l=\infty$}
The next result is an upper bound on the mixing time at 
$\l=\infty$, showing that the estimate of Theorem \ref{unilobo} is
sharp up to constant factors in this case. 
\begin{Proposition}
\la{lainf}
For $\l=\infty$, both with and without the wall
\be
\tmix\leq 
\,L^2\,.
\la{lobo14}
\end{equation}
\end{Proposition}
\proof
We first give the proof for the system with the wall.  Let $A(\h)$
denote the area under the path $\eta$:
\begin{eqnarray}
\label{eq:area}
A(\eta):=\sum_{x\in\Lambda}\eta_x.
\end{eqnarray}
Then, $A(\eta^\wedge(t))$ is a process on
$\{A(\vee),\dots,A(\wedge)\}$, where $A(\wedge)=\frac{L^2}4$ is the
maximal value and $A(\vee)=\frac{L}2$ is the minimal value. The
process starts at $A(\wedge)$, has $\pm 2$ increments and is killed
upon hitting $A(\vee)$. We want an upper bound on the expected value
of $\t$, where $\t$ denotes the hitting time of $A(\vee)$. It will be
shown below that
\be
\cL A(\eta) \leq -1\,,\quad \forall \eta\neq \vee\,. 
\la{lobo15}
\end{equation}
Assume (\ref{lobo15}) and consider the martingale
\be
M_t =
A(\eta_t) - A(\wedge) - \int_0^t\cL A(\eta_s) ds\,,
\la{lobo015}
\end{equation}
where $\eta_t:=\eta^\wedge(t)$. 
By the optional stopping theorem and (\ref{lobo15}) 
we obtain
\be
0 = \bbE M_\t =  \bbE A(\eta_\t) - A(\wedge) - \bbE \int_0^\t\cL A(\eta_s) ds
\geq  A(\vee) - A(\wedge) + \bbE\t\,.
\la{lobo155}
\end{equation}
This implies $\bbE\t \leq A(\wedge)-A(\vee)$ 
and therefore, using (\ref{varcoup})
and Markov's inequality:
\begin{align*}
\|P_t(\wedge,\cdot) - \d_{\vee}\|_{\rm var} &\leq 
\bbP(\t>t) \nonumber\\
&\leq \frac1t\,\left(A(\wedge)-A(\vee)\right)\,.
\end{align*}
This gives the mixing time bound 
$\tmix\leq \nep{}\,(A(\wedge)-A(\vee))\leq 
\,\frac{\nep{}}4\,
L^2\leq L^2$. 

It remains to prove (\ref{lobo15}). From Lemma \ref{lemma1}, with
$\d=0$, 
we have 
\begin{align}
\cL A(\eta)
&= 
 \sum_x \left(\Delta \eta\right)_x 
+ \sum_{x=1}^{L-1}\left[1_{\{\eta_{x-1}=\eta_{x+1}=0\}}-
1_{\{\eta_{x-1}=\eta_{x+1}=1\}}\right]\,.
\la{lobo17}
\end{align}
Note that for any $\eta\in\O_L^+$ we have $\sum_x \left(\Delta
  \eta\right)_x = -1$ (for a non-negative path the number of mountains exceeds by 1
the number of valleys, deterministically). 
The last term in (\ref{lobo15}) can be
estimated by observing that whenever $\eta\neq \vee$ then the number
of sites $x$ such that  $\eta_{x-1}=\eta_{x+1}=0$ is at most equal to the
number of sites $x$ such that  $\eta_{x-1}=\eta_{x+1}=1$. It follows
that  
$$
\sum_{x=1}^{L-1}\left[
1_{\{\eta_{x-1}=\eta_{x+1}=0\}}-
1_{\{\eta_{x-1}=\eta_{x+1}=1\}}\right] \leq 1_{\{\eta = \vee\}}
\,.
$$ 
In particular, 
$$
\cL A(\eta)\leq - \, 1_{\{\eta \neq \vee\}}\,.
$$
This ends the proof of (\ref{lobo15}).

\medskip

Finally, we prove the proposition for the system without the wall.
Here the equilibrium measure 
$\mu$ is the uniform probability on all $2^{\frac{L}2}$
configurations 
$$
\O_0=\{\eta\in\O_L\,:\;\;\eta_x=0 \;\;\text{for all even}\; 
x\}\,.$$ Let $\t$ denote the hitting time of $\O_0$ for our process
$\eta$ started in $\wedge$. Since $\l=\infty$ the process cannot exit
$\O_0$ once it has entered. It is then obvious that $\t$ coincides
with the first time when the configuration started from $\wedge$ and
evolving {\sl with the wall} and $\l=\infty$ equals $\vee$, the
(zigzag) minimal configuration satisfying the hard--wall
constraint. Therefore, we already know from (\ref{lobo015}) that
$\bbE\,\t \leq L^2/4$. Let $-$ and $+$ denote the minimal and maximal
configurations in $\O_0$ respectively, and write $\eta^-(t),\eta^+(t)$
for the associated evolutions. The standard coupon--collector estimate
gives that the coupling time $\t'$ of $\eta^-(t),\eta^+(t)$ satisfies
$\bbE\t' = O(\log L)$ (there are $L/2$ independent coordinates to be
updated). Let now $\wt \t$ denote the coupling time for
$\eta^\wedge(t)$ and $\eta^-(t)$. We see that
\begin{align*}
\max_{\eta\in\O_L}
\|P_t(\eta,\cdot) - \mu\|_{\rm var}
& \leq \|P_t(\wedge,\cdot) - P_t(-,\cdot)\|_{\rm var}\\
& \leq\bbP(\wt \t\geq t) \leq \frac1t\, \bbE[\wt \t] \\
& \leq \frac1t\, (\bbE[\t] + \bbE[\t']) 
\leq  \frac1t\, (L^2/4 + O(\log L))\,.
\end{align*}
As before, this gives the mixing time bound 
$\tmix\leq L^2$, provided $L$ is sufficiently large.

\section{Proof of Theorem \ref{th5}}
\la{52}
\bigskip
The statement concerning $\l\geq1$ has been proven 
in Section \ref{sec:gap1L} (spectral gap upper bound) and in Section 
\ref{sec:universal} (mixing time lower bound), so we only need to prove
\eqref{th5_2}. Before we do that, we give a heuristic argument which 
suggests that the $L^{-5/2}$ behavior in Theorem \ref{th5} might be
the correct one.

\subsection{A heuristic justification of the $L^{-5/2}$ result}
\label{sec:miclo}
Consider the model without wall and $\l<1$, and start the dynamics
from a non-negative initial configuration $\xi$, e.g.,
$\xi=\wedge$. We know that the equilibrium measure $\mu^\l_L$ is
symmetric under $\h\leftrightarrow -\h$ and, from \eqref{eq:nozeri}, that
$\mu^\l_L(\eta\leq0)>0$ uniformly in $L$. Also, from the analysis of the
model with the wall, we know that the dynamics restricted to
configurations $\h\geq 0$ (or to $\h\leq0$) relaxes in a time of order
at most $O(L^2\log L)$. Therefore, it is reasonable that the
relaxation time of our system {\sl without wall} is of the same order as
the first time $\tau$ such that $\h^\xi_x(\tau)\leq0$ for every $x$,
{\sl provided that} $\tau\gg L^2\log L$. On the other hand, it is
plausible that the most convenient mechanism for the system to go from
an initial configuration $\xi\geq0$ to some $\eta\leq 0$ is the
following:
\begin{enumerate}
\item first of all  a ``negative bubble'' is formed 
close to one of the borders of the system (say, the left border). By
``negative bubble at the left border'' we mean that there exists
$0<x\leq L/2$ such that $\h_y\leq 0$ for $y\leq 2x$ and $\h_y\geq 0$
for $y\geq2x$. The point $2x$ will be referred to as the {\sl right-hand
boundary of the bubble}.  Of course, when the bubble is first created
one has $2x=2$
\item the bubble grows until it occupies the whole $\L$, \ie
until $2x=L$.
\end{enumerate}
Processes involving several bubbles or the formation of a bubble far
away from the system boundaries would require that the configuration
$\eta^\xi$ develops more zeros, and therefore they look much less
likely in view of Eq. \eqref{eq:pochizeri} and \eqref{eq:farfrom0L};
at any rate, we neglect them. Now we introduce a simplified model
which mimics the process of bubble formation and growth described
above. We will implicitly assume that at any time $t$ the system is at
equilibrium conditionally on the position, $2x$, of the right--hand
border of the bubble. Again, this is reasonable provided that $\tau\gg
L^2\log L$.
\begin{figure}[h]
  \centerline{
  \psfrag{x}{$2x$}
  \psfrag{0}{$0$}
  \psfrag{L}{$L$}
\psfig{file=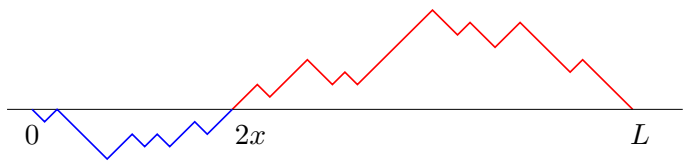,height=0.8in}}
  \caption{A typical configuration with a negative bubble
  at the left border. Apart from the 
  point $2x$ (the right-hand border of the bubble) the polymer has very few
  zeros, since the line is repulsive.}
 \label{fig:bolla}
 \end{figure}

Consider a birth-death process on $\{0,\ldots,L\}$ with invariant measure
\begin{eqnarray}
\label{eq:nu}
\nu(x):=Z^{-1} \frac1{(x\vee 1)^{3/2}((L-x)\vee 1)^{3/2}},
\end{eqnarray}
where $Z=Z(L)$ normalizes $\nu(\cdot)$ to $1$. It is clear that 
$Z\approx L^{-3/2}$ and 
\begin{eqnarray}
\label{eq:nuapprox}
\nu(x)=\nu(L-x)\approx (x\vee 1)^{-3/2}\;\;\mbox{if}\;\;x\leq L/2,
\end{eqnarray}
where $A\approx B$ means that there exists a universal constant $c$
such that $(1/c)\leq A/B\leq c$.  We consider a Metropolis dynamics where
 the ``birth'' rate, $b(x)$, of jump from $x$ to $x+1$ is given
for $x<L$ by $\min(1,\nu(x+1)/\nu(x))$, while the death rates are
uniquely determined by the requirement that $\nu(\cdot)$ be
reversible. 

The connection of this dynamics with the ``bubble dynamics''
discussed above is obvious if one interprets $2x$ as the rightmost point
of a bubble in a system of length $2L$, in view of
\begin{eqnarray}
\mu^\l_{2L}(\h_y\leq 0\;\mbox{for}\; 
y<2x;\h_y\geq0\;\mbox{for}\;y>2x)=\frac{Z^+_{2x}(\l)Z^+_{2(L-x)}(\l)}
{Z_{2L}(\l)}
\end{eqnarray}
and of Theorem \ref{th:asymptZ}.

The following two observations will be useful in a while:
\begin{eqnarray}
\label{eq:osserv1}
\sum_{y= x}^L\nu(y)\approx 1\approx\sum_{y=0}^ x\nu(y)
\end{eqnarray}
and
\begin{eqnarray}
\label{eq:osserv2}
b(x)\approx 1\,.
\end{eqnarray}
We will estimate how the inverse spectral gap of the birth-death
process, $\gap(L)^{-1}$, grows with $L$ applying a method of Hardy
inequalities due to L. Miclo \cite{cf:Miclo}. For this we need some
additional notation, and we define for $0\leq i\leq L$
\begin{eqnarray}
\label{eq:B+}
B_+(i)&:=&\sup_{x>i}\left(\sum_{y=i+1}^x\frac1{\nu(y)b(y)}\right)\sum_{y\geq
x}\nu(y)
\\
\label{eq:B-}
B_-(i)&:=&\sup_{x<i}\left(\sum_{y=x}^{i-1}\frac1{\nu(y)b(y)}\right)
\sum_{y\leq x}\nu(y)
\\
\label{eq:B}
B&:=&\min_{0\leq i\leq L}\left(B_+(i)\vee B_-(i)\right),
\end{eqnarray}
with the convention that $B_+(L)=B_-(0)=0$.
Then, Proposition 3.1 of \cite{cf:Miclo} says that
\begin{eqnarray}
\frac B2\leq \gap(L)^{-1}\leq 4B.
\end{eqnarray}
In view of \eqref{eq:osserv1} and \eqref{eq:osserv2}, if we are only
interested in the order of magnitude of the inverse spectral gap as a
function of $L$ and not in precise constants, we can replace $b(y)$,
$\sum_{y\geq x}\nu(y)$ and $\sum_{y\leq x}\nu(y)$ by $1$ in
\eqref{eq:B+} and  \eqref{eq:B-}.  Using
\eqref{eq:nuapprox}, one finds
\begin{eqnarray}
B_+(i)\approx B_-(L-i)\approx
\left\{
\begin{array}{lll}
L^{5/2}&\mbox{if}& i\leq L/2\\
(L-i+1)^{5/2}&\mbox{if}& i\geq L/2
\end{array}
\right.,
\end{eqnarray}
which immediately implies that $\gap(L)^{-1}\approx B\approx L^{5/2}$.
Note that, in contrast with 
Theorem \ref{th5}, no spurious logarithmic factor appears.
Note also that the equilibration time for this birth-death process is 
indeed much larger than $L^2\log L$, as required for the heuristic 
argument to be consistent, see discussion before Eq. \eqref{eq:nu}.

\subsection{Proof of bound \eqref{th5_2}}

We need some preliminary notation.  Let
  $w:\bbR\ni x\mapsto w(x)\in[0,1]$ be a smooth function such that
  $w(x)=1$ for $x\leq -1$ and $w(x)=0$ for $x>1$. Recall the definition 
\eqref{eq:area} of $A(\eta)$.
\begin{Theorem}
\label{th:5/2} Let $\l<1$ 
and define
  \begin{eqnarray}
    f(\h):=w\left(\frac{A(\h)}{L^{3/2}(\log L)^{-3}}\right).
  \end{eqnarray}
There exists $c(\l,w)<\infty$ such that 
\begin{eqnarray}
\label{eq:5/2}
\frac{ {\mathcal E}(f,f)}{\var(f)}
< c(\l,w)\frac{(\log L)^8}{L^{5/2}}.
\end{eqnarray}
\end{Theorem}
As a consequence, we will deduce that if we start from the maximal
configuration $\wedge$ then at any given time $t\ll L^{5/2}/(\log L)^8$ the
area $A(\eta^\wedge(t))$ is larger than $L^{3/2}/(\log L)^{3}$ with large
probability.  More precisely:
\begin{Proposition}
\label{th:prop52} For every $\l<1$ 
  \begin{eqnarray}
\label{eq:prop52}
\liminftwo{L\to\infty\,,}{t=o(L^{5/2}(\log L)^{-8})}
\bbP\left(A(\eta^\wedge(t))\geq L^{3/2}(\log L)^{-3}\right)=1.
  \end{eqnarray}
\end{Proposition}
Of course, since at equilibrium $ \mu^\l_L(A\leq0)\geq1/2$, this
implies directly that the mixing time in this situation is at least  
$\Omega(L^{5/2}(\log L)^{-8})$.

\smallskip

\noindent
{\sl Proof of Theorem \ref{th:5/2}} For notational simplicity let
$\e_L:=(\log L)^{-3}$. We put also 
$$
b_L:=(\log L)^{-8/3}
$$
and
$$
\gamma_L:=(\log L)^{-7/6}.
$$ 

It is easy to show that the variance of $f$ converges to $1/4$
 for $L\to\infty$. Indeed, for $L$ large the function
$w(A(\eta)/(L^{3/2}\e_L))$ takes the value $1$ with
$\mu^\l_L$--probability $1/2+o(1)$ and the value $0$ also with
probability $1/2+o(1)$. This is quite intuitive from the properties of
the delocalized phase discussed in section \ref{sec:delocph}, but more
precisely it follows from
\begin{eqnarray}
1/2&\geq& \mu^\l_L\left( A(\eta)\geq L^{3/2}\e_L\right)
= \mu^\l_L\left( A(\h)\leq -L^{3/2}\e_L\right) \\\nonumber
&\geq& \frac12
\mu_L^\l\left(\nexists x:\;L\, b_L<x<L-L\, b_L,\eta_x=0\right)\\\nonumber
&&\times\mu_L^\l\left(\left. |A(\h)|\geq L^{3/2}\e_L\right|
\nexists x:\;L\, b_L<x<L-L\, b_L,\eta_x=0\right)
\end{eqnarray}
together with  Eq. \eqref{eq:farfrom0L} and the fact that the last factor
in the right-hand side is bounded below by $(1-2/L)$ for large $L$, 
 (cf.\ Lemma \ref{th:lemma52} below and the subsequent discussion).
 Therefore,
\begin{eqnarray}
\var(f)=1/4+o(1).
\end{eqnarray}
As for the Dirichlet form,
\begin{eqnarray}
 {\mathcal E}(f,f)\leq 4\frac{\left(\max_{x\in[-1,1]}|w'(x)|
 \right)^2}
{L^{2}\e_L^2}
\mu^\l_L\left(\left|A(\h)\right|\leq L^{3/2}\e_L\right).
\end{eqnarray}
A factor $L$ comes from the sum over $x$ in \eqref{eq:dirich}, while
the factor $L^{-3}/{\e_L^2}$ originates from 
$$
\left|f(\h)-Q_x f(\eta)\right|
\le
\frac{2}{\e_L\,L^{3/2}}\times \max_{x\in\left[\frac{A(\h)-1}{L^{3/2}\e_L},
\frac{A(\h)+1}{L^{3/2}\e_L}
\right]}|w'(x)|.
$$
In order to conclude the proof, it is therefore sufficient to prove 
that
\begin{eqnarray}
\label{eq:A32}
\mu^\l_L\left(\left|A(\h)\right|\leq L^{3/2}\e_L\right)=O(L^{-1/2}(\log L)^2).
\end{eqnarray}
To this end, observe first of all that 
\begin{eqnarray}
\label{eq:decompo}
\mu^\l_L\left(\left|A(\h)\right|\leq L^{3/2}\e_L\right)
&\leq& \mu^\l_L\left(\exists x:\;L\, b_L<x<L-L\, b_L,\eta_x=0\right)
\\\nonumber
&+&\mu^\l_L\left(\left|A(\h)\right|\leq L^{3/2}\e_L;
\nexists x:\;L\, b_L<x<L-L\, b_L,\eta_x=0\right).
\end{eqnarray}
The first term in the right-hand side of \eqref{eq:decompo} is of
order $O(L^{-1/2}(\log L)^{4/3})$, thanks to \eqref{eq:farfrom0L} and to the
definition of $b_L$. As for the second one, decompose for convenience
\begin{align*}
  A(\h)&=\sum_{1\leq x<L\,b_L}\eta_x+\sum_{L\,b_L\leq x\leq L-L\,b_L}\eta_x+
\sum_{L\,b_L<x<L}\eta_x\\
&=:A^{(1)}(\h)+A^{(2)}(\h)+A^{(3)}(\h)\,.
\end{align*}
The key estimate we need is 
\begin{Lemma}
\label{th:lemma52}
For $L$ sufficiently large
one has
\begin{eqnarray}
\label{eq:rho1}
     \mu^\l_L\left(|A^{(1)}(\h)|+|A^{(3)}(\h)|\geq (1/2)L^{3/2}\e_L\right)
\leq \frac1L
  \end{eqnarray}
and
\begin{eqnarray}
\label{eq:rho2}
    \mu^\l_L\left(|A^{(2)}(\h)|\leq 2L^{3/2}\e_L;\nexists
   x:\;L\,b_L<x<L-L\,b_L,\eta_x=0\right)\leq \frac1L.
\end{eqnarray}
\end{Lemma}
Indeed, thanks to the lemma we obtain immediately from \eqref{eq:decompo}
\begin{eqnarray}
\mu^\l_L\left(\left|A(\h)\right|\leq L^{3/2}\e_L\right)
\leq O(L^{-1/2}(\log L)^2)+2/L
\end{eqnarray}
which concludes the proof of the theorem. \qed


\smallskip

\noindent
{\sl Proof of Lemma \ref{th:lemma52}}. 
The proof of \eqref{eq:rho1} is easy. We start by observing that
\begin{eqnarray*}
 \mu^\l_L\left(|A^{(1)}(\h)|\geq (1/4)L^{3/2}\e_L\right)
&\leq&\mu^\l_L\left (\max_{x<L\,b_L}|\eta_x|\geq (1/4)L^{1/2}\frac{\e_L}
{b_L}
\right)
\\\nonumber
&=&\mu^{+,2\l}_L\left(\max_{x<L\,b_L}\eta_x\geq (1/4)L^{1/2}\frac{\e_L}
{b_L}\right)
\\\nonumber&\leq&
\mu^{+,0}_L\left(\max_{x<L\,b_L}\eta_x\geq (1/4)L^{1/2}\frac{\e_L}
{b_L}\right)
\end{eqnarray*}
where we used monotonicity (say, FKG) in the last inequality.  Since
$\mu^{+,0}_L(\cdot)=\bP(\cdot\tc\eta_L=0,\;\eta_x>0\;\forall\; 1<x<L)$
where we recall that  $\bP(\cdot)$ the law of the one-dimensional 
simple random walk started at $0$, we have
\begin{eqnarray}
 \mu^\l_L\left(|A^{(1)}(\h)|\geq (1/4)L^{3/2}\e_L\right)
\leq \frac{
\bP(\max_{x<L\,b_L}\eta_x\geq (1/4)L^{1/2}\e_L/b_L)}
{\bP(\eta_L=0;\;\eta_x>0\;\forall\; 1<x<L)}.
\end{eqnarray}
  Now, for the denominator we employ
\eqref{eq:denomina}, while for the numerator we observe that
\begin{align*}
&\bP\left(\max_{x<L\,b_L}\eta_x\geq (1/4)L^{1/2}\e_L/b_L\right)
=\sum_{\ell\geq  (1/4)L^{1/2}\e_L/b_L}
\bP\left(\max_{x<L\,b_L}\eta_x=\ell\right)
\\
&\quad\quad\quad\quad\quad\quad = \sum_{\ell\geq  (1/4)L^{1/2}\e_L/b_L}
\max(\bP\left(\eta_{L\,b_L}=\ell\right),
\bP\left(\eta_{L\,b_L}=\ell+1\right)),
\end{align*}
where we used  \cite[Sec. III.7,Th. 1]{cf:Feller} in the last equality.
From this one sees that
\begin{eqnarray}
 \mu^\l_L\left(|A^{(1)}(\h)|\geq (1/4)L^{3/2}\e_L\right)\leq
\frac1{2L}
\end{eqnarray}
for $L$ large.
Since $A^{(1)}(\h)$ and $A^{(3)}(\h)$ are equally distributed, this proves
\eqref{eq:rho1}. 

As for \eqref{eq:rho2}, we note (using also the symmetry
$\eta\leftrightarrow -\eta$) that the left-hand side is bounded above
by
\begin{eqnarray}
\mu^\l_L\left(0<A^{(2)}(\h)\leq 2 L^{3/2}\e_L| \eta_x>0\;\forall
\;L\,b_L<x<L-L\,b_L\right),
\end{eqnarray}
which by FKG is itself bounded above for $L$ large by
\begin{eqnarray}
\label{eq:areapiccola}
\bP\left(A_\ell(\eta)\leq 4 \ell^{3/2}\e_\ell\tc\eta_\ell=0;\eta_x>0 \;
\forall\; 1<x<\ell\right)
\end{eqnarray}
where we put $\ell:=\ell(L):=L-2\lfloor L\,b_L\rfloor$ and for clarity
of notation $A_\ell(\eta):=\sum_{x\leq\ell}\eta_x$.  Letting
$M_\ell:=\lceil\ell\,\gamma_\ell\rceil$, one can bound above
\eqref{eq:areapiccola} by
\begin{eqnarray}
\label{eq:Meta}
&&\bP\left(\max_{x\leq \ell}\eta_x\leq \sqrt{M_\ell}
\tc\eta_\ell=0;\eta_x>0 \;\forall\; 1<x<\ell\right)\\\nonumber&&+
\bP\left(\max_{x\leq \ell}\eta_x> \sqrt{M_\ell}
;A_\ell(\eta)\leq 4 \ell^{3/2}\e_\ell
\tc\eta_\ell=0;\,\eta_x>0 \;\forall\; 1<x<\ell\right).
\end{eqnarray}
Using monotonicity twice, the first term of
\eqref{eq:Meta} is easily bounded above by
\begin{eqnarray}
\nonumber
&&\bP\left(\max_{x\leq \ell}\eta_x\leq \sqrt{M_\ell}
\tc\eta_1>0,\ldots,\eta_{\ell-1}>0;\eta_\ell=0;
\eta_{2j M_\ell}=2,\,j\leq\lfloor
1/(4\gamma_\ell)\rfloor
\right)\\\nonumber
&&\leq\left[\bP\left(\eta_{ M_\ell}\leq 
\sqrt{M_\ell}\tc\eta_1>0,\ldots,\eta_{2M_\ell-1}>0;\eta_{2
M_\ell}=0\right)\right]^{1/(4\gamma_\ell)}
\\\nonumber
&&\leq
\left[\bP\left(\eta_{M_\ell}\leq 
 \sqrt{M_\ell}\tc\eta_{2M_\ell}=0\right)\right]^{1/(4\gamma_\ell)}
\end{eqnarray}
Since $n^{-1/2}\eta_{n/2}$ converges weakly for $n\to\infty$ under
$\bP(\cdot\tc\eta_{n=0})$ to a non-degenerate Gaussian random variable (the
Brownian Bridge at time $1/2$), the probability in the last expression
is strictly smaller than $1$ uniformly in $\ell$, and therefore the
first term in \eqref{eq:Meta} is smaller than $1/L$ for $L$ large.

As for the second term in \eqref{eq:Meta}, we note that the conditions
on $\max_{x\leq \ell}\eta_x$ and on $A_\ell(\h)$ imply that there exist
$1<x,y<\ell$ such that $|x-y|\leq4\ell\e_\ell/\sqrt{\gamma_\ell}$
and $|\eta_x-\eta_y|\geq (1/2)\sqrt{M_\ell}$ (just 
take as $x$ the position of the maximum of $\eta$). As a consequence,
using \eqref{eq:denomina} one can bound above the second term
in \eqref{eq:Meta} by
\begin{eqnarray}
c\,\ell^{3/2}\bP\left(\exists 1<x<y<\ell: |x-y|\leq
4\ell\e_\ell/\sqrt{\gamma_\ell}, |\eta_x-\eta_y|\geq
(1/2)\ell^{1/2}\sqrt{\gamma_\ell}\right)
\end{eqnarray}
and for this quantity the upper bound $1/\ell$
 for $\ell$ large follows immediately from
standard simple-random-walk estimates. The factor $\ell^{3/2}$ arises from
the estimate \eqref{eq:denomina}. \qed

\bigskip

\noindent
{\sl Proof of Proposition \ref{th:prop52}.}  
Let $w_+:\bbR\ni
x\mapsto w_+(x)\in[0,1]$ be a smooth function such that $w_+(x)=0$
for $x\leq 1/2$ and $w_+(x)=1$ for $x\geq1$, and define $w_-(.)$ via
$w_-(x)=w_+(-x)$. We put $f_\pm(\eta):=w_\pm(A(\eta)/L^{3/2}\e_L)$
where, as in the proof of Theorem \ref{th:5/2}, $\e_L:=(\log L)^{-3}$.
The proof of Theorem \ref{th:5/2} can be repeated essentially without changes
to show that 
\begin{eqnarray}
\label{eq:nochanges}
 {\mathcal E}(f_\pm,f_\pm)=O(L^{-5/2}(\log L)^8).
\end{eqnarray}

To begin the proof of \eqref{eq:prop52}, observe that by monotonicity
\begin{eqnarray}
 \bbP \left(A(\eta^\wedge(t))\geq L^{3/2}\e_L\right)
\geq \frac{\int d\mu^\l_L(\xi)\, \left(f_+(\xi)\right)^2
 \bbP \left(A(\eta^\xi(t))\geq L^{3/2}\e_L\right)
}
{\mu^\l_L\left( (f_+)^2\right)
}.
\end{eqnarray}
It is immediate to realize that 
\begin{eqnarray}
\label{eq:immediate}
 \mu^\l_L\left((f_+)^2\right)=
\mu^\l_L\left(f_+\right)+o(1)=1/2+o(1)
\end{eqnarray}
for $L\to\infty$. Using reversibility of
dynamics and Cauchy-Schwarz in the numerator one obtains then
\begin{eqnarray}
\nonumber
 \bbP \left(A(\eta^\wedge(t))\geq L^{3/2}\e_L\right)&\geq& (2+o(1))
\int d \mu^\l_L(\xi)\,{\bf 1}_{\{A(\xi)\geq L^{3/2}\e_L\}}\,
\left[\left( P_t f_+\right)(\xi)\right]^2
\\
\label{eq:spezzata}
&=&(2+o(1))
\mu^\l_L\left[
\left( P_t f_+\right)^2\right]
\\\nonumber
&-&(2+o(1))
\int d\mu^\l_L(\xi)\,{\bf 1}_{\{A(\xi)< L^{3/2}\e_L\}}\,
\left[\left( P_t f_+\right)(\xi)\right]^2.
\end{eqnarray}
We will show later that
\begin{eqnarray}
\label{eq:tchebi}
\var\left ( P_t f_+\right)\geq \var
(f_+)e^{-2t\frac{c(\l)(\log L)^8}{L^{5/2}}}.
\end{eqnarray}
From \eqref{eq:immediate} one then deduces that 
\begin{eqnarray}
(2+o(1))\mu^\l_L\left[
\left( P_t f_+\right)^2\right]\geq 1+o(1)
\end{eqnarray}
for $t=o(L^{5/2}/(\log L)^8)$. 
As for the integral in \eqref{eq:spezzata}, rewrite it as
\begin{eqnarray}
\label{eq:integral}
\int d\mu^\l_L(\xi)\,{\bf 1}_{\{|A(\xi)|< L^{3/2}\e_L\}}\,
\left[\left( P_t f_+\right)(\xi)\right]^2+
\int d\mu^\l_L(\xi)\,{\bf 1}_{\{A(\xi)\leq - L^{3/2}\e_L\}}\,
\left[\left( P_t f_+\right)(\xi)\right]^2.
\end{eqnarray}
The first term is $o(1)$ as follows from \eqref{eq:A32} plus the 
fact that $f_+$ is bounded. The second one, on the other hand, is 
bounded above by
\begin{eqnarray}
\la{eq:scalprod}
\mu^\l_L\left(\,f_-\,
 P_t f_+\right)
\end{eqnarray}
(indeed, recall that $||f_+||\leq 1$.)
It is obvious from the definition of $f_\pm$ that this integral vanishes
at $t=0$. To show that \eqref{eq:scalprod} is $o(1)$ we evaluate the
$t$--derivative of it: using reversibility and  Cauchy-Schwarz,
\begin{eqnarray}
\left|\frac{d}{dt}\mu^\l_L\left(f_-\,
 P_t f_+\right)\right|&=&
\left|\mu^\l_L\left( f_-\,(-{\mathcal L})
P_t f_+\right)\right|\\\nonumber
&=& \left|\mu^\l_L\left((-{\mathcal L})^{1/2}f_-\,
(-{\mathcal L})^{1/2}P_t f_+\right)\right|\\\nonumber 
&\leq& \sqrt{
 {\mathcal E}(f_-,f_-)\, {\mathcal E}(f_+,f_+)
}=O(L^{-5/2}(\log L)^8).
\end{eqnarray}
We can therefore conclude that \eqref{eq:integral} is $o(1)$ for $
t=o(L^{5/2}/(\log L)^8)$.

Finally, we prove \eqref{eq:tchebi}. This is a simple consequence of the
general inequality
\begin{eqnarray}
\var\left(P_t f\right)\geq \var(f) e^{-2t\frac{\mathcal E(f,f)}
{\var(f)}},
\end{eqnarray}
which holds for every $f$ thanks to the spectral theorem plus Jensen's
inequality, and of Eqs. \eqref{eq:nochanges} and \eqref{eq:immediate}. \qed


\bigskip

\section{Further results in the localized phase}\la{locco}
In this section we prove Theorem \ref{th:relaxV} and Theorem
\ref{th:t13}. All our arguments below refer to the system with the
wall with $\l>2$.

%

\subsection{Proof of Theorem \ref{th:relaxV}} \label{sec:relaxV}
Recall the definition (\ref{eq:Boltz_wall}) of the equilibrium measure
$\mu^{+,\l}_L$ and set
\begin{eqnarray}
U(L,t):=\max_{0<x<L}\left(\mu^{+,\l}_L(\eta_x)-
\mathbb E_L(\eta^\vee_x(t)\right)\geq 0\,,
\end{eqnarray}
where for later convenience we indicated explicitly the
$L$--dependence in the average over the process. Non-negativity
follows from monotonicity. Also, from monotonicity and Markov's inequality
we have 
\be
\|P_t(\vee,\cdot)-\mu^{+,\l}_L\|_{\rm var}\leq \frac12\,L\,U(L,t)\,.
\la{tovar}
\end{equation}

Let $\ell=\ell(L):=2 \lfloor c_0\log L\rfloor\in 2\bbN$ where $c_0$ will be 
chosen sufficiently large later.
Thanks 
to the exponential decay of correlations (cf.\ Lemma
\ref{th:tight} and subsequent discussion), 
one has for every $\ell/2\leq x\leq L-\ell/2$
\begin{eqnarray}
\label{eq:expc1}
0\leq\mu^{+,\l}_{L}(\eta_x)- \mu^{+,\l}_{\ell}(\eta_{\ell/2})\leq c\, e^{-
\ell/c}\,.
\end{eqnarray}
Here and below we write $c$ for a suitable constant, whose value may
vary from line to line.
For $1\leq x\leq \ell/2$ one has  instead
\begin{eqnarray}
\label{eq:expc2}
0\leq\mu^{+,\l}_{L}(\eta_x)- \mu^{+,\l}_{\ell}(\eta_{x})\leq c\,
e^{-\ell/c}\,,
\end{eqnarray}
and for $L-\ell/2\leq x\leq L$
\begin{eqnarray}
\label{eq:expc3}
0\leq\mu^{+,\l}_{L}(\eta_x)- \mu^{+,\l}_{\ell}(\eta_{x-L+\ell})\leq c\,
e^{-\ell/c}\,.
\end{eqnarray}
If \eg 
$\ell/2\leq x\leq L-\ell/2$, then 
\eqref{eq:expc1} implies 
\begin{eqnarray}
0&\leq& \mu^{+,\l}_{L}(\eta_x)-\bbE_{L}(\eta^\vee_x(t))\leq c\,
e^{-\ell/c}+\mu^{+,\l}_\ell(\eta_{\ell/2})-\bbE_{L} (\eta^\vee_x(t))\\\nonumber
&\leq&c\,e^{-\ell/c}+\mu^{+,\l}_\ell(\eta_{\ell/2})-
\bbE_{\ell} (\eta^\vee_{\ell/2}(t))\,,
\end{eqnarray}
where we used again monotonicity in the last inequality.
For $x\not\in [\ell/2,L-\ell/2]$ one obtains analogous bounds from
Eqs. \eqref{eq:expc2}--\eqref{eq:expc3}.
As a consequence, one concludes that for every $t\geq0$
\begin{eqnarray}\la{ecco1}
U(L,t)\leq c\, e^{-\ell(L)/c}+U(\ell(L),t)\,.
\end{eqnarray}
From \eqref{claims} it follows
that for every $n\in 2\bbN$
\begin{eqnarray*}
|\mu^{+,\l}_n(\eta_x)-\bbE_n(\eta_x^\xi(t))|\leq c\, n^3 e^{-t/(c\, n^2)}\,,
\end{eqnarray*}
for every $x\leq n,t\geq0$ and every initial condition $\xi$. 
Therefore, (\ref{ecco1}) implies
\begin{eqnarray*}
U(L,t)\leq \frac{c}{L^{c_0/c}}+c\, c_0^3 (\log L)^3 e^{-t/(c\, c_0^2\, 
(\log L)^2)}\,.
\end{eqnarray*}
If $t\geq t_0(L):=c_1(\log L)^3$ with $c_1$ sufficiently large, then 
\begin{eqnarray}
\label{eq:Ult}
U(L,t)\leq \frac{c}{L^{c_0/c}}\,.
\end{eqnarray}
Finally, if we choose $c_0$ sufficiently large in the definition 
of $\ell(L)$, it follows from \eqref{eq:Ult} and (\ref{tovar}) that 
for $t\geq t_0(L)$ the variation distance between $\mu^{+,\l}_L(.)$ and the 
distribution of $\eta^\vee(t)$ is $o(1)$, and Eq. \eqref{eq:UB_V} is proven.

\smallskip

Let us now turn to the proof of the lower bound on the
equilibration time starting from $\vee$. It is possible to apply the
ideas of
\cite{cf:hayes_sinclair} to prove that the dynamics starting from
$\vee$ takes at least a time of order $\log L$ to relax to equilibrium,
but we shall prove the stronger statement \eqref{eq:LB_V}. To begin, we define 
$\mathcal C$ to be the
set $\mathcal C:=\{2j\lfloor \sqrt L\rfloor, j\leq \lfloor \sqrt L\rfloor/2
-1\}$ and 
$$ 
f(\eta):=\frac1{|\mathcal C|}\sum_{x\in \mathcal C}\,1_{\eta_x=0}\,.  
$$ 
Using
the exponential decay of correlations (Lemma \ref{th:tight})
we see that
\begin{eqnarray}
\label{eq:varf}
\var(f)\leq \frac {c}{\sqrt L}\,,
\end{eqnarray}
where the variance is computed w.r.t.\ $\mu^{+,\l}_L$. 
Next, we need the following estimate, whose proof will be given later. 
\begin{Lemma}
\label{th:lemmatlogt}
There exist positive constants $c_0=c_0(\l)$ and $L_0$ such
that for every $x\in \mathcal C$, $0<t<\sqrt L$ and $L>L_0$ one has
\begin{eqnarray}
\label{eq:mediestaccate}
\bbP( \eta_x^\vee(t)=0)- \mu^{+,\l}_L(\h_x=0)
\geq c_0\, e^{-\sqrt{t}/c_0}.
\end{eqnarray}
\end{Lemma}
In order to prove that at time $(\log L)^2/c$ the total variation
distance from equilibrium is still $1+o(1)$, we introduce the set
$$
K:=\{\eta\in\Omega^+_L:\;f(\eta)\geq \mu^{+,\l}_L(f)+(c_0/2)\, e^{- \sqrt
{t}/c_0}\},
$$ 
where $c_0$ is the same as in \eqref{eq:mediestaccate},
and we  show that $\mu^{+,\l}_L(K)=o(1)$ while 
\begin{eqnarray}
\label{eq:hnotinK}
\bbP(\eta^\vee(t)\in K)=1+o(1). 
\end{eqnarray}
The first fact follows from \eqref{eq:varf} and
Chebyshev's inequality:
\begin{eqnarray}
\mu^{+,\l}_L(K)\leq \frac{c}{c_0^2\sqrt L}e^{2\sqrt {t}/c_0}=o(1)\,,
\end{eqnarray}
as $L\to\infty$, if $t\leq (\log L)^2/c$.  As for \eqref{eq:hnotinK}, it is convenient
to introduce a modified process, call it $\tilde \h^\vee(t)$, which
is just the original process started  from $\vee$ but 
conditioned on the event  that, for every $t\geq0$, $\h^\vee_x(t)\equiv
\vee_x$ for every $x$ such that $\min\{|x-j|:j\in \mathcal C\}\geq
\sqrt L/2$ (in other words, the points at
distance at least $\sqrt L/2$ from $\mathcal C$ are kept at their
initial values for all times). Denote by $\pi_t(\cdot)$ (respectively
$\tilde\pi_t(\cdot)$) the marginal distribution of
$\{\eta_x\}_{x\in\mathcal C}$ under the law of $\eta^\vee(t)$
(resp. the law of $\tilde \h^\vee(t)$). The proof of the next claim is
postponed for a moment.
\begin{claim}
\label{th:claimpipitilde}
\begin{eqnarray}
\label{eq:pipitilde}
||\pi_t(\cdot)-\tilde\pi_t(\cdot)||_{\rm var}\leq c\,e^{-(\log L)\,\sqrt L/c}
\end{eqnarray}
if $ t\leq (\log L)^2/c$.
\end{claim}
The same is actually true as long as $t\leq L^{1/2-\e}$, but we will
not need that. Assuming the validity of (\ref{eq:pipitilde}) we are
able to finish the proof of the theorem.

Observe that 
\begin{eqnarray}
\label{eq:varpitilde}
\var_{\tilde\pi_t}(f)\leq \frac c{\sqrt L}\,,
\end{eqnarray}
since, from the way the modified dynamics is constructed,
 $\tilde\h_x^\vee(t)$ is independent of $\tilde\h_y^\vee(t)$ for
$x,y \in \mathcal C$ with $x\ne y$.
From Eqs. \eqref{eq:pipitilde} and \eqref{eq:varpitilde} and the fact that
$||f||_\infty= 1$ one deduces 
that
\begin{eqnarray}\la{ecco2}
\var_{\pi_t}(f)\leq \frac c{\sqrt L}\,.
\end{eqnarray}
Thanks to (\ref{ecco2}) and \eqref{eq:mediestaccate}, equation 
\eqref{eq:hnotinK} is seen to hold for all $t\leq (\log
L)^2/c$ by an application of Chebyshev's inequality.
\qed

\smallskip

\noindent
{\sl Proof of Claim \ref{th:claimpipitilde}}. This is based on a
standard disagreement percolation argument, see {\sl e.g.} 
\cite[Sec. 3.1]{cf:hayes_sinclair}. 
Consider $n$ IID Poisson clocks of rate $1$ and let $p(n,t)$ be the
probability that there is an increasing sequence of times
$0<t_1<\ldots<t_n<t$ such that the clock labeled $i$ rings at time
$t_i$. An standard computation gives that
\begin{eqnarray}
p(n,t)<\left(\frac{et}n\right)^n.
\end{eqnarray}
On the other hand, it is immediate to realize that
\begin{eqnarray}
\bbP\left(\exists x\in \mathcal C:\; \h_x^\vee(t)\ne \tilde \h_x^\vee(t)
\right)\leq 2|\mathcal C|\,p(\lfloor \sqrt L/2\rfloor,t)\,,
\end{eqnarray}
from which Eq. \eqref{eq:pipitilde} easily follows.
\qed

\smallskip

\begin{figure}[h]
\centerline{
\psfrag{x}{$x$}
\psfrag{xl}{$x+\ell$}
\psfrag{xl-}{$x-\ell$}
\psfig{file=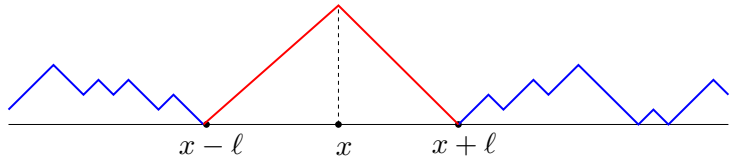,height=0.8in}}
\caption{A typical path belonging to $B_{x,\ell}$.}
\label{fig:Bxt}
\end{figure}

\noindent
{\sl Proof of Lemma \ref{th:lemmatlogt}}.
Define, for $\ell\in2\bbN$, the set
\begin{eqnarray}
\label{eq:Qt}
B_{x, \ell}:=
\{\xi: \xi_{x+j}=\ell-|j|\mbox{\;for every\;}
|j|\leq \ell\}.
\end{eqnarray}
In other words, configurations $\xi\in B_{x,\ell}$ take in
$\{x-\ell,\ldots,x+\ell\}$ the maximal value allowed by the
constraint $\xi_{x\pm\ell}=0$ (see Fig. \ref{fig:Bxt}).

Let $\ell(t)\in2\bbN$ satisfy
\begin{eqnarray}
\label{eq:cepsilonbis}
c_1 \sqrt{t}<\ell(t)<2c_1\sqrt{t}\,,
\end{eqnarray}
for some sufficiently large constant $c_1$ to be chosen later.
Since we are in the
localized phase,
\begin{eqnarray}
\label{eq:pKt}
\mu^{+,\l}_L(B_{x,\ell(t)})\geq c\,e^{-\sqrt{t}/c}\,,
\end{eqnarray}
uniformly in $L$.
 Indeed, it is not difficult to deduce from Lemma \ref{th:tight}
that $\mu^{+,\l}_L(\h_a=\h_b=0)$ is bounded away from zero uniformly in 
$a,b,L\in 2\bbZ$. On the other hand, from the definition of the model
\begin{eqnarray*}
\mu^{+,\l}_L(B_{x,\ell(t)})=\frac{\mu^{+,\l}_L(\h_{x\pm\ell(t)=0})}
{Z^+_{2\ell(t)}(\l)},
\end{eqnarray*}
from which the claim \eqref{eq:pKt} follows through Theorem \ref{th:asymptZ}.

\smallskip

Next, we write
\begin{eqnarray}
\label{eq:possottint}
\bbP( \eta_x^\vee(t)=0)- \mu^{+,\l}_L(\h_x=0)
=\int d\mu_L^{+,\l}(\xi)\,\left[\bbP( \eta_x^\vee(t)=0)
-\bbP(\eta^\xi_x(t)=0)\right]\,.
\end{eqnarray}
Since the quantity which is being integrated in the right-hand side of
\eqref{eq:possottint} is non-negative by monotonicity, 
(\ref{eq:pKt}) implies 
\begin{align*}
&\bbP( \eta_x^\vee(t)=0)- \mu^{+,\l}_L(\h_x=0)\\
& \quad\geq c\,e^{-\sqrt t/c}\,
\int d\mu_L^{+,\l}(\xi\tc \xi\in B_{x,\ell(t)})\,
\left[\bbP( \eta_x^\vee(t)=0)
-\bbP(\eta^\xi_x(t)=0)\right]\,.
\end{align*}
Note that, if $\xi\in
B_{x,\ell(t)}$ then $\eta^\xi(s)$ is stochastically higher, for every
$s>0$, than the configuration $\hat\h^\xi(s)$ which has law
$\bbP(\cdot|\eta^\xi_{x\pm\ell(t)}(r)=0\;\forall \;r\leq t)$.
This holds in particular for $s=t$.
Therefore, 
\begin{eqnarray}
\bbP( \eta_x^\vee(t)=0)- \mu^{+,\l}_L(\h_x=0)\geq 
 c\,e^{-\sqrt t/c}\,
\left[\mu_L^{+,\l}(\eta_x=0)
-\bbP \left( \hat\eta^\xi_x(t)=0\right)
\right]\,,
\end{eqnarray}
where it is clear that the last term is independent of the choice of
$\xi\in B_{x,\ell(t)}$. Indeed, $\{\hat \h^\xi_x(s)\}_{s\geq0}$ depends only
on the value of $\xi$ in the interval $\{x-\ell(t),x+\ell(t)\}$, on
which however there is no choice once we require that $\xi\in
B_{x,\ell(t)}$.  

Next, we shall use the following estimate, the proof of which is
postponed for a moment. Recall that $c_1$ is the constant defining
$\ell(t)$ in \eqref{eq:cepsilonbis}.
\begin{claim}
\label{clamsi}
For any $\e_1,\e_2 >0$, there exists $C>0$ such that for all $c_1\geq C$:
\be
\label{clamsi1}
\bbP(\hat\eta^\xi_x(t) < (1-\e_1)\ell(t))\leq \e_2\,,\quad t\geq 0\,.
\end{equation}
\end{claim}

From (\ref{clamsi1}), for any given $\e>0$, if $c_1$
is chosen sufficiently large as a function of
$\e$, we have $$
\bbP\left( \hat\eta^\xi_x(t)=0\right)<\e\,.
$$ 
Choosing
$\e<\mu_L^{+,\l}(\eta_x=0)/2$, the desired estimate
\eqref{eq:mediestaccate} follows. Note that this $\e>0$ can be chosen
to be independent of $L$ since in the localized phase the
probabilities
$ \mu_L^{+,\l}(\eta_x=0)$ are uniformly bounded away from zero.
This ends the proof of the Lemma \ref{th:lemmatlogt}.
\qed

\smallskip

\noindent
{\sl Proof of Claim \ref{clamsi}}.
By monotonicity it is sufficient to prove the claim at $\l=\infty$. 
Let $\varphi_L(t)$ denote the height in the middle of the segment 
$\{0,\dots,L\}$, at time $t$, of the usual process $\eta^\wedge(t)$.
It suffices to prove that, for every $\e_1,\e_2>0$ there exists $\d>0$ such that
\be\la{eso2}
\bbP\left(\varphi_L(t)\leq (1-\e_1)\frac{L}2\right)
\leq \e_2\,,\quad\;t\leq \d\,L^2\,.
\end{equation}
This can be shown to follow from the argument in the proof of Theorem
\ref{unilobo}. Namely, let $B_D$ and $\xi(t)$ be the area and the
auxiliary process defined there. Note that the geometric construction
of Section \ref{sec:universal} implies, 
in particular, that if $\varphi_L(t)\leq (1-\e_1)\frac{L}2$ then 
$B_D(\xi(t))\leq (1-\e')\,B_D(\wedge)$ for some $\e'>0$.
Also, recall that $\bbE B_D(\xi(t))\geq B_D(\wedge) - t$, so that 
$$
\var[B_D(\xi(t))]=\bbE [B_D(\xi(t))^2] - \bbE [B_D(\xi(t))]^2\leq B_D(\wedge)^2 -
[B_D(\wedge)-t]^2\leq 2tB_D(\wedge)\,.
$$
Then the claim follows from an application of Chebyshev's inequality.
\qed

\subsection{Stretched exponential decay of local observables:
Proof of Theorem \ref{th:t13}}
\label{sec:locfunct}
We start with the proof of the upper bound (\ref{eq:t13}). 
We first prove it under the extra assumption that 
$f$ is monotone decreasing and non-negative.  
For $\ell_1,\ell_2\in2\bbZ$ with $
\ell_1<\ell_2$ let $$\mu^+_{\ell_1,\ell_2}(\cdot):=
\mu^+_\infty(\cdot|\eta_{\ell_1}=\h_{\ell_2}=0)$$ 
denote the equilibrium measure with zero boundary conditions at
$\ell_1$ and $\ell_2$. From Lemma
\ref{th:tight} it follows that if we choose $\ell_i$ such that $\mathcal
S_f\subset
\{\ell_1,\ldots,\ell_2\}$, then
\begin{eqnarray}
  \left|\mu^+_\infty(f)-\mu^+_{\ell_1,\ell_2}(f)\right|\leq c\,
  ||f||_\infty \left[ e^{-d(\mathcal S_f,\{\ell_1\})/c}+e^{-d(\mathcal
  S_f,\{\ell_2\})/c}\right].
\end{eqnarray}
By positivity and monotonicity of 
$f$ one has for every $t\geq 0$ that
\begin{eqnarray*}
  \mu^+_\infty
((P_t f)^2)\leq \mu^+_{\ell_1,\ell_2}((P_t f)^2)
\leq \int d\mu^+_{\ell_1,\ell_2}(\xi)\left(\left. \bbE\left(f(\h^\xi(t))\right|
\h^\xi_{\ell_1}(r)=\h^\xi_{\ell_2}(r)=0\;\forall\;r\leq t\right)
\right)^2\,.
\end{eqnarray*}
Therefore, using the lower bound on the
spectral gap given by 
\eqref{th1_1}, 
one finds
\begin{eqnarray}
\nonumber
  \var_{\mu^+_\infty}(P_t f)&\leq&
  \var_{\mu^+_{\ell_1,\ell_2}}(f)e^{-t/[c(\ell_2-\ell_1)^2]}
\\
&&+ c\, ||f||^2_\infty \left[
 e^{-d(\mathcal S_f,\{\ell_1\})/c}+e^{-d(\mathcal
  S_f,\{\ell_2\})/c}\right]\,.\la{ecco45}
\end{eqnarray}
We may
choose 
$\ell_i, i=1,2$ such that 
$$t^{1/3}\leq
d(\mathcal S_f,\{\ell_i\})< 2 t^{1/3}\,,$$ 
which gives $(\ell_2-\ell_1)<c t^{1/3}$ for $t^{1/3}>\diam(\mathcal S_f)$.
Since $\var_{\mu^+_{\ell_1,\ell_2}}(f)\leq \|f\|_\infty^2$,
(\ref{ecco45}) proves
\eqref{eq:t13} for all 
$t$ such that $t>(\diam(\mathcal S_f))^3$.  
If $t$ is smaller than that then we obtain again the claimed bound by
adjusting
the constant $C_f$. This proves \eqref{eq:t13} for $f$ bounded, local,
non-negative and decreasing. 

To prove the claim for any bounded local $f$ we first introduce a
cutoff
parameter $\ell_0$ and rewrite
$f$ as $f=f_0 + f_1$ where $f_0(\eta)=f(\eta)1_E$ and
$f_1(\eta)=f(\eta)1_{E^c}$
with $E$ representing the event $\{\max_{x\in
    \cS_f}\eta_x \leq \ell_0\}$. Observe that 
\be
\var_{\mu^+_\infty}(P_t f)\leq 2 \var_{\mu^+_\infty}(P_t f_0) +
2\|f\|_\infty^2\,
\mu^+_\infty(E^c)\,.
\la{ecco6}
\end{equation}
Since in the localized phase the height at any point has an
exponential tail, one has $\mu^+_\infty(E^c)\leq
c\,|\cS_f|\,\nep{-\ell_0/c}$, where $|\cS_f|$ stands for the
cardinality of $\cS_f$. Let now $\O_0$ denote the set of all possible 
values of the configuration $\{\eta_x, x\in\cS_f\}$ that are compatible
with the constraint $\eta\in E$. Note that its cardinality $|\O_0|$ 
is at most $C\,\ell_0$, for some constant $C$ depending on $\cS_f$.
We can write $1_{\{\eta_x=\si_x\}} = 1_{\{\eta_x\leq\si_x\}} -
1_{\{\eta_x<\si_x\}}$ and expand
$$
f_0(\eta) = \sum_{\si\in\O_0} f(\si)
\,\prod_{x\in\cS_f}1_{\{\eta_x=\si_x\}}
= \sum_{\si\in\O_0}\sum_{A\subset\cS_f} (-1)^{|\cS_f\setminus A|}\,f(\si) \,g_{\si,A}\,,
$$ 
where $g_{\si,A}:= \prod_{x\in A,y\in\cS_f\setminus
  A}1_{\{\eta_x\leq\si_x\}}1_{\{\eta_y<\si_y\}}$. The latter is a
bounded local, non-negative and decreasing function to which the
argument leading to (\ref{ecco45}) applies.  
Adjusting the constant $C_f$ we may therefore
estimate
\begin{align*}
\var_{\mu^+_\infty}(P_t f_0)& \leq
C_f\,\ell_0\,\sum_{\si\in\O_0}\sum_{A\subset\cS_f}
\var_{\mu^+_\infty}(P_t g_{\si,A})\\
&\leq C_f\,\ell_0\,\nep{-t^{1/3}/c}\,.
\end{align*}
Recalling (\ref{ecco6}), it suffices to take $\ell_0=t^{1/3}$ to
conclude the proof.

\smallskip

We turn to the proof of the lower bound (\ref{eq:t12}).  
Let $f=f^{\underline a,I}$ be a function as in (\ref{eq:cylinder}). 
Assume that
\begin{eqnarray}
\label{eq:tgrande12}
\sqrt{t}>2\min_{x\in I}a_x.
\end{eqnarray}
and let $y\in I$ be a point such that
$a_y=\min_{x\in I}a_x$ (to fix ideas, we assume that $y$ is even). Let 
$\ell(t)\in2\bbN$ satisfy (\ref{eq:cepsilonbis})
for some sufficiently large $c_1$.
We write
\begin{eqnarray}
\var_{\mu^+_\infty}(P_t f)&=&\frac12 \mu^{+}_\infty\otimes \mu^+_\infty
\left[((P_t f)(\xi)-(P_t f)
(\xi'))^2\right]\\\nonumber
&\geq& \frac12\mu^+_\infty\otimes\mu^+_\infty
\left[1_{\{\xi\in B_{y,\ell(t)}\}}\left((P_t f)(\xi)-
(P_t f)(\xi')\right)^2
\right]
\end{eqnarray}
where $B_{y,\ell}$ is the set defined in \eqref{eq:Qt}.
As a consequence of \eqref{eq:pKt},
\begin{eqnarray}
\label{eq:Kt}
\var_{\mu^+_\infty}(P_t f)\geq c\,e^{-\sqrt{t}/c}\mu^+_\infty{\otimes}
\mu^+_\infty
\left[\left.
\left((P_t f)(\xi)-(P_t f)(\xi')\right)^2
\right|\xi\in B_{y,\ell(t)}\right].
\end{eqnarray}
By monotonicity, for every initial condition $\xi\in B_{y,\ell(t)}$ and every
$s>0$ (and in particular for $s=t$) one has
\begin{eqnarray}
0\leq(P_t f)(\xi)\leq \bbE\left[\left. f( \h^\xi(s))\right|
\h^\xi_{y\pm\ell(t)}(r)=0\;\forall \;r\leq s\right].
\end{eqnarray}
From (\ref{clamsi1}) we know that
for any
given $\e>0$, if $c_1=c_1(\e)$ is chosen large enough in
\eqref{eq:cepsilonbis}, for every $s\leq t$ one has
\begin{eqnarray}
\bbP\left[\left.\eta^\xi_y(s)\geq \ell(t)/2\right|
\h^\xi_{y\pm\ell(t)}(r)=0\;\forall \;r\leq s
\right]\geq 1-\e
\end{eqnarray}
if $\xi\in B_{y,\ell(t)}$.  Since $||f||_\infty=1$ and $\ell(t)>2 a_y$
(cf.\ \eqref{eq:tgrande12}), this implies that if $\xi\in
B_{y,\ell(t)}$, then 
$$ 0\leq(P_t f)(\xi)\leq
\e.  $$ 
Going back to \eqref{eq:Kt}, we obtain
\begin{eqnarray}
\var_{\mu^+_\infty}(P_t f)\geq c\,e^{-\sqrt{t}/c}\left[
\mu^+_\infty\left[(P_t f)^2\right]-2\e\right] \,.
\end{eqnarray}
Choosing $\e$ small enough and using the Cauchy-Schwarz inequality we find the
estimate \eqref{eq:t12}. Adjusting the value of $c_f$ yields
the desired bound for all $t\geq 0$, \ie without the restriction
\eqref{eq:tgrande12}. 
This ends the proof of Theorem \ref{th:t13}.
\qed

\smallskip

{\bf Acknowledgments} F.L.T.\ would like to thank 
C.\ Bernardin, T.\ Bodineau, F.\ Caravenna and G.\ Giacomin 
for useful discussions
and suggestions on the literature. P.C.\ thanks N.\
Yoshida, P.\ Tetali and D.\ Randall for useful conversations. 
F.L.T.\ was partially supported by the GIP-ANR project JC05\_42461
(POLINTBIO).

\end{document}